\documentclass[a4paper,11pt]{article}
\usepackage{graphicx}
\usepackage{amsfonts}
\usepackage{amsmath}
\usepackage{amssymb}
\usepackage{indentfirst}
\usepackage{titlesec}

\def\barr{\begin{array}}
\def\earr{\end{array}}
\def\bali{\begin{aligned}}
\def\eali{\end{aligned}}
\def\bearr{\begin{eqnarray}}
\def\eearr{\end{eqnarray}}

\providecommand{\play}{\displaystyle}

\providecommand{\li}{\limits}

\providecommand{\pt}{\partial}

\providecommand{\ra}{\rightarrow}

\providecommand{\da}{\downarrow}

\providecommand{\Prob}{\mathbf P}

\providecommand{\E}{\mathbf E}

\providecommand{\al}{\alpha}

\providecommand{\bt}{\beta}

\providecommand{\gm}{\gamma}

\providecommand{\Gm}{\Gamma}

\providecommand{\dt}{\delta}

\providecommand{\Dt}{\Delta}

\providecommand{\ve}{\varepsilon}

\providecommand{\kp}{\kappa}

\providecommand{\lb}{\lambda}

\providecommand{\Lb}{\Lambda}

\providecommand{\sm}{\sigma}

\providecommand{\om}{\omega}

\providecommand{\Om}{\Omega}

\providecommand{\R}{\mathbb R}

\providecommand{\Z}{\mathbb Z}

\providecommand{\cE}{\mathcal E}

\providecommand{\cF}{\mathcal F}

\providecommand{\1}{\mathbf 1}

\providecommand{\grad}{\nabla}

\providecommand{\torus}{\mathbb{T}}

\providecommand{\boldq}{\boldsymbol{q}}

\providecommand{\boldp}{\boldsymbol{p}}

\providecommand{\boldW}{\boldsymbol{W}}

\providecommand{\boldb}{\boldsymbol{b}}

\providecommand{\boldal}{\boldsymbol{\al}}

\providecommand{\boldbt}{\boldsymbol{\bt}}

\providecommand{\boldgm}{\boldsymbol{\gm}}

\providecommand{\boldR}{\boldsymbol{R}}

\providecommand{\boldx}{\boldsymbol{x}}

\providecommand{\boldX}{\boldsymbol{X}}

\providecommand{\boldy}{\boldsymbol{y}}

\providecommand{\boldN}{\boldsymbol{N}}

\providecommand{\bolde}{\boldsymbol{e}}

\providecommand{\boldphi}{\boldsymbol{\varphi}}

\begin{document}

\author{Mark Freidlin\thanks{Dept of Mathematics, University of Maryland at
College Park, mif@math.umd.edu.}, Wenqing Hu\thanks{Dept of
Mathematics, University of Maryland at College Park,
huwenqing@math.umd.edu.}}
\title{Smoluchowski-Kramers approximation in the case of variable friction}
\date{}

\maketitle
\bigskip

\begin{abstract}
We consider the small mass asymptotics (Smoluchowski-Kramers
approximation) for the Langevin equation with a variable friction
coefficient. The limit of the solution in the classical sense does
not exist in this case. We study a modification of the
Smoluchowski-Kramers approximation. Some applications of the
Smoluchowski-Kramers approximation to problems with fast oscillating
or discontinuous coefficients are considered.
\end{abstract}



\section{Introduction}

The Langevin equation $$\mu
\ddot{\boldq}_t^\mu=\boldb(\boldq_t^\mu)-\lb
\dot{\boldq}_t^\mu+\sm(\boldq_t^\mu)\dot{\boldW}_t \ , \
\boldq_0^\mu=\boldq\in \R^n \ , \dot{\boldq}_0^\mu=\boldp\in \R^n \
, \eqno(1.1)$$ describes the motion of a particle of mass $\mu$ in a
force field $\boldb(\boldq)$, $\boldq\in \R^n$, subjected to random
fluctuations and to a friction proportional to the velocity. Here
$\boldW_t$ is the standard Wiener process in $\R^n$, $\lb>0$ is the
friction coefficient. The vector field $\boldb(\boldq)$ and the
matrix function $\sm(\boldq)$ are assumed to be continuously
differentiable and bounded together with their first derivatives.
The matrix $a(\boldq)=(a_{ij}(\boldq))=\sm(\boldq)\sm^*(\boldq)$ is
assumed to be non-degenerate.

Put $\boldp_t^\mu=\dot{\boldq}_t^\mu$. Then (1.1) can be written as
a first order system:

$$\left\{\begin{array}{l}
\play{\dot{\boldq}_t^\mu=\boldp_t^\mu}
\\
\play{\dot{\boldp}_t^\mu=\dfrac{1}{\mu}\boldb(\boldq_t^\mu)
-\dfrac{\lb}{\mu}\boldp_t^\mu+\dfrac{1}{\mu}\sm(\boldq_t^\mu)\dot{\boldW}_t}
\end{array}\right. \eqno(1.2)$$

The diffusion process $(\boldp_t^\mu, \boldq_t^\mu)=\boldX_t^\mu$ in
$\R^{2n}$ is governed by the generator $L$:

$$Lu(\boldp,\boldq)=\dfrac{1}{2\mu^2}\sum\li_{i,j=1}^n a_{ij}(\boldq)\dfrac{\pt^2 u}{\pt p_i\pt p_j}
+\dfrac{1}{\mu}(\boldb(\boldq)-\lb\boldp)\cdot
\grad_{\boldp}u+\boldp \cdot \grad_{\boldq}u \ .$$

Note that, since functions $\boldq_t^\mu$ are continuously
differentiable with probability one,

$$\int_0^t \sm_{ij}(\boldq_s^\mu)dW_s^j=\sm_{ij}(\boldq_t^\mu)
W_t^j-\int_0^t W_s^j\left(\grad_{\boldq}\sm_{ij}(\boldq_s^\mu) \cdot
\boldp_s^\mu\right)ds \ .$$

This allows to consider equations (1.2) for each trajectory
$\boldW_t$ individually, and there is no necessity in the
introduction of a stochastic integral. In particular, if (1.2) is
considered as a stochastic differential equation, stochastic
integrals in the It\^{o} and in the Stratonovich sense coincide:
$\play{\int_0^t \sm(\boldq_s^\mu)d\boldW_s=\int_0^t
\sm(\boldq_s^\mu)\circ d\boldW_s}$.

It is assumed usually that the friction coefficient $\lb$ is
constant. Under this assumption, one can prove that $\boldq_t^\mu$
converges in probability as $\mu \da 0$ uniformly on each finite
time interval $[0,T]$ to an $n$-dimensional diffusion process
$\boldq_t$: for any $\kp, T>0$ and any $\boldp_0^\mu=\boldp\in \R^n$
fixed,
$$\lim\li_{\mu \da 0}\Prob\left(\max\li_{0\leq t \leq T}
|\boldq_t^\mu-\boldq_t|_{\R^d}>\kp\right)=0 \ .$$

Here $\boldq_t$ is the solution of equation

$$\dot{\boldq}_t=\dfrac{1}{\lb}\boldb(\boldq_t)+
\dfrac{1}{\lb}\sm(\boldq_t)\dot{\boldW}_t  \ , \
\boldq_0=\boldq_0^\mu=\boldq\in \R^n \ . \eqno(1.3)$$

The stochastic term in (1.3) should be understood in the It\^{o}
sense.

The approximation of $\boldq_t^\mu$ by $\boldq_t$ for $0<\mu<<1$ is
called the Smoluchowski-Kramers approximation. This is the main
justification for replacement of the second order equation (1.1) by
the first order equation (1.3). The price for such a simplification,
in particular, consists of certain non-universality of equation
(1.3): The white noise in (1.1) is an idealization of a more regular
stochastic process $\dot{\boldW}_t^\dt$ with correlation radius
$\dt<<1$ converging to $\dot{\boldW}_t$ as $\dt \da 0$. Let
$\boldq_t^{\mu,\dt}$ be the solution of equation (1.1) with
$\dot{\boldW}_t$ replaced by $\dot{\boldW}_t^\dt$. Then limit of
$\boldq_t^{\mu,\dt}$ as $\mu,\dt \da 0$ depends on the relation
between $\mu$ and $\dt$. Say, if first $\dt \da 0$ and then $\mu \da
0$, the stochastic integral in (1.3) should be understood in the
It\^{o} sense; if first $\mu \da 0$ and then $\dt \da 0$,
$\boldq_t^{\mu,\dt}$ converges to the solution of (1.3) with
stochastic integral in the Stratonovich sense. (See, for instance,
[5].)

Consider now the case of a variable friction coefficient
$\lb=\lb(\boldq)$. We assume that $\lb(\boldq)$ has continuous
bounded derivatives and $0<\lb_0\leq \lb(\boldq)\leq \Lb<\infty$. It
turns out, as we will see in the next section, that in this case the
solution $\boldq_t^\mu$ of (1.1) does not converge, in general, to
the solution of (1.3) with $\lb=\lb(\boldq)$, so that the
Smoluchowski-Kramers approximation should be modified. In order to
do this, we consider equation (1.1) with $\dot{\boldW_t}$ replaced
by $\dot{\boldW}_t^\dt$ described above:

$$\mu \ddot{\boldq}_t^{\mu,\dt}=\boldb(\boldq_t^{\mu,\dt})-
\lb(\boldq_t^{\mu,\dt})\dot{\boldq}_t^{\mu,\dt}
+\sm(\boldq_t^{\mu,\dt})\dot{\boldW}_t^\dt \ , \
\boldq_0^{\mu,\dt}=\boldq \ , \ \dot{\boldq}_0^{\mu,\dt}=\boldp \ .
\eqno(1.4)$$

We prove that after such a regularization, the solution of (1.4) has
a limit $\widetilde{\boldq}_t^\dt$ as $\mu \da 0$, and
$\widetilde{\boldq}_t^\dt$ is the unique solution of the equation
obtained from (1.4) as $\mu=0$:

$$\dot{\widetilde{\boldq}}_t^\dt=
\dfrac{1}{\lb(\widetilde{\boldq}_t^\dt)}\boldb(\widetilde{\boldq}_t^\dt)
+
\dfrac{1}{\lb(\widetilde{\boldq}_t^\dt)}\sm(\widetilde{\boldq}_t^\dt)\dot{\boldW}_t^\dt
 \ , \ \widetilde{\boldq}_0^\dt=\boldq \ . \eqno(1.5)$$

Now we can take $\dt \da 0$ in (1.5). As the result we get the
equation

$$\dot{\widehat{\boldq}}_t=\dfrac{1}{\lb(\widehat{\boldq}_t)}\boldb(\widehat{\boldq}_t)+
\dfrac{1}{\lb(\widehat{\boldq}_t)}\sm(\widehat{\boldq}_t)\circ
\dot{\boldW}_t \ , \ \widehat{\boldq}_0=\boldq \ , \eqno(1.6)$$
where the stochastic term should be understood in the Stratonovich
sense. So the regularization leads to a modified
Smoluchowski-Kramers equation (1.6). We prove this in Section 3.

Some applications of the Smoluchowski-Kramers approximation are
considered in Sections 4 and 5: the case of fast oscillating in the
space variable, periodic or stochastic, friction coefficient is
studied; gluing condition at the discontinuity points of the
friction coefficient are considered. In the last Section 6, we
briefly consider some remarks and generalizations.

\

\textbf{Notations.} We use $|\bullet|_{\R^d}$ to denote the standard
Euclidean norm in $\R^d$. When $d=1$ we set
$|\bullet|_{\R^1}=|\bullet|$. For a vector-valued function
$\boldsymbol{f}(\boldx)=(f_1(\boldx),...,f_d(\boldx))$, $\boldx\in
\R^d$, we set $\|\boldsymbol{f}\|_{\infty}=\max\li_{1\leq i\leq
d}\|f_i\|_{\infty}=\max\li_{1\leq i\leq d}\sup\li_{\boldx\in
\R^d}|f_i(\boldx)|$. All the vectors are marked with either bold
letters or with an arrow on it.

\

\section{Some estimates. The classical Smoluchowski-Kramers
approximation does not work for variable friction coefficients}

We consider the following system

$$\mu \ddot{\boldq}^\mu_t=\boldb(\boldq_t^\mu)-\lb(\boldq_t^\mu)\dot{\boldq}_t^\mu
+\dot{\boldW}_t \ , \ \boldq_0^\mu=\boldq\in \R^d
 \ , \ \dot{\boldq}_0^\mu=\boldp\in \R^d. \eqno(2.1)$$
Here $\infty>\Lb\geq\lb(\bullet)\geq \lb_0 >0$ is a function of
$\boldq_t^\mu$. We assume that function $\lb(\bullet)$ and the
vector field $\boldb(\bullet)$ are continuously differentiable and
bounded together with their first derivatives. The process
$\boldW_t$ is the standard Wiener process in $\R^d$. For simplicity
of calculations we consider here the case when the diffusion matrix
$a(\bullet)$ is the identity (compare with (1.1)). The case of
general diffusion matrix can be considered in a similar way and we
will briefly mention it in Section 6.

Let $\boldp_t^\mu=\dot{\boldq}_t^\mu$, we have, that (2.1) is
equivalent to the system

$$
\left\{\begin{array}{l} \dot{\boldq}_t^\mu=\boldp_t^\mu  \ , \\
\dot{\boldp}_t^\mu=\dfrac{1}{\mu}\boldb(\boldq_t^\mu)-
\dfrac{\lb(\boldq_t^\mu)}{\mu}\boldp_t^\mu+\dfrac{1}{\mu}\dot{\boldW}_t
\ .
\end{array}
\right.\eqno(2.2)$$

Then

$$\dfrac{d}{dt}\left(e^{\frac{1}{\mu}\int_0^t
\lb(\boldq_s^\mu)ds}\boldp_t^\mu\right)=e^{\frac{1}{\mu}\int_0^t
\lb(\boldq_s^\mu)ds}\left(\dot{\boldp}_t^\mu+\dfrac{1}{\mu}\lb(\boldq_t^\mu)\boldp_t^\mu\right)
=e^{\frac{1}{\mu}\int_0^t
\lb(\boldq_s^\mu)ds}\left(\dfrac{1}{\mu}\boldb(\boldq_t^\mu)+\dfrac{1}{\mu}\dot{\boldW}_t\right)
\ ,$$ and

$$e^{\frac{1}{\mu}\int_0^t \lb(\boldq_s^\mu)ds}\boldp_t^\mu - \boldp
=\dfrac{1}{\mu}\int_0^t e^{\frac{1}{\mu}\int_0^s
\lb(\boldq^\mu_r)dr} \boldb(\boldq_s^\mu) ds + \frac{1}{\mu}\int_0^t
e^{\frac{1}{\mu}\int_0^s \lb(\boldq_r^\mu)dr}d\boldW_s \ .
\eqno(2.3)$$

For notational convenience we introduce the function
$\play{A(\mu,t)=\int_0^t \lb (\boldq_s^\mu) ds}$. It is clear that
$t\Lb\geq A(\mu,t)\geq t \lb_0$. Using (2.3) we have, that

$$\boldp_t^\mu=e^{-\frac{1}{\mu}A(\mu,t)}\left(\boldp+ \dfrac{1}{\mu}\int_0^t e^{\frac{1}{\mu}A(\mu,s)}\boldb(\boldq_s^\mu)ds
+ \frac{1}{\mu}\int_0^t e^{\frac{1}{\mu}A(\mu,s)}d\boldW_s\right)\ .
$$

Therefore we have

$$\begin{array}{l}
\play{\boldq_t^\mu=\boldq+\int_0^t \boldp_s^\mu ds}
\\
\play{=\boldq+\boldp\int_0^t e^{-\frac{1}{\mu}A(\mu,s)}ds+
\dfrac{1}{\mu}\int_0^t e^{-\frac{1}{\mu}A(\mu,s)}\left(\int_0^s
e^{\frac{1}{\mu}A(\mu,r)}\boldb(\boldq_r^\mu)dr\right)ds+}
\\
\play{ \ \ \ \ \ \ \ \ \ \ + \dfrac{1}{\mu}\int_0^t
e^{-\frac{1}{\mu}A(\mu,s)}\left(\int_0^s
e^{\frac{1}{\mu}A(\mu,r)}d\boldW_r\right)ds }
\\
\play{ = \boldq + \boldal(\mu) + \boldbt(\mu) + \boldgm(\mu)\
.}\end{array} \eqno(2.4)$$ Here $\boldal(\mu), \boldbt(\mu),
\boldgm(\mu)$ are three (vector) functions in the right hand side of
(2.4):

$$\boldal(\mu)=\boldp\int_0^t e^{-\frac{1}{\mu}A(\mu,s)}ds \ ,$$

$$\boldbt(\mu)=\dfrac{1}{\mu}\int_0^t e^{-\frac{1}{\mu}A(\mu,s)}\left(\int_0^s
e^{\frac{1}{\mu}A(\mu,r)}\boldb(\boldq_r^\mu)dr\right)ds \ ,$$

$$\boldgm(\mu)=\dfrac{1}{\mu}\int_0^t
e^{-\frac{1}{\mu}A(\mu,s)}\left(\int_0^s
e^{\frac{1}{\mu}A(\mu,r)}d\boldW_r\right)ds \ .$$

In the following we will use the relation
$$\dfrac{d}{dt}\left(e^{-\frac{1}{\mu}A(\mu,t)}\right)
=-\dfrac{1}{\mu}e^{-\frac{1}{\mu}A(\mu,t)}\dfrac{dA(\mu,t)}{dt}
=-\dfrac{1}{\mu}e^{-\frac{1}{\mu}A(\mu,t)}\lb(\boldq_t^\mu) \ .
\eqno(2.5)$$

We will also use the estimates
$$\dfrac{\mu}{c\Lb}(1-e^{-\frac{c\Lb t}{\mu}})=\int_0^t e^{-\frac{c\Lb
s}{\mu}}ds\leq\int_0^t e^{-\frac{c}{\mu}A(\mu,s)}ds\leq \int_0^t
e^{-\frac{c\lb_0 s}{\mu}}ds=\dfrac{\mu}{c\lb_0}(1-e^{-\frac{c\lb_0
t}{\mu}})\leq\dfrac{\mu}{c\lb_0} \ , \eqno(2.6)$$

$$\begin{array}{l}
\play{\dfrac{\mu}{c\Lb}(1-e^{-\frac{c\Lb t}{\mu}})=\int_0^t
e^{-\frac{c\Lb (t-s)}{\mu}}ds\leq \int_0^t
e^{-\frac{c}{\mu}(A(\mu,t)-A(\mu,s))}ds\leq}
\\
\play{ \ \ \ \ \ \ \ \ \ \ \ \ \ \ \ \leq \int_0^t e^{-\frac{c\lb_0
(t-s)}{\mu}}ds=\dfrac{\mu}{c\lb_0}(1-e^{-\frac{c\lb_0 t}{\mu}})\leq
\dfrac{\mu}{c\lb_0} \ .}
\end{array} \eqno(2.7)$$

Here $c$ is a positive constant.

\

We get in this section some bounds for $\boldal(\mu)$,
$\boldbt(\mu)$, $\boldgm(\mu)$ which show, in particular, that the
classical Smoluchowski-Kramers approximation does not hold in the
case of variable friction. These bounds also will be used to obtain
a modified Smoluchowski-Kramers approximation.

\

\textbf{2.1. Estimates of $\boldal(\mu)$. }

\

We have, by (2.5),

$$\begin{array}{l}
\play{\boldal(\mu)=\boldp\int_0^t e^{-\frac{1}{\mu}A(\mu,s)}ds}
\\
\play{=\boldp\int_0^t(-\mu)\dfrac{1}{\lb(\boldq_s^\mu)}d(e^{-\frac{1}{\mu}A(\mu,s)})}
\\
\play{=-\boldp\mu\left[\dfrac{e^{-\frac{1}{\mu}A(\mu,t)}}{\lb(\boldq_t^\mu)}-\dfrac{1}{\lb(\boldq)}-\int_0^t
e^{-\frac{1}{\mu}A(\mu,s)}d(\frac{1}{\lb(\boldq_s^\mu)})\right]} \ .
\end{array}$$

Let
$R_{\boldal}(\mu)=\mu\left[\dfrac{e^{-\frac{1}{\mu}A(\mu,t)}}{\lb(\boldq_t^\mu)}-\dfrac{1}{\lb(\boldq)}\right]$.
It is easy to estimate $|R_{\boldal}(\mu)|\leq \dfrac{\mu}{\lb_0}$.
Therefore $|R_{\boldal}(\mu)|\ra 0$ as $\mu \da 0$.

Let $$(I)=\int_0^t
e^{-\frac{1}{\mu}A(\mu,s)}d(\frac{1}{\lb(\boldq_s^\mu)}) . $$ We
have

$$\begin{array}{l}
\play{(I)=-\int_0^t
e^{-\frac{1}{\mu}A(\mu,s)}\dfrac{1}{\lb^2(q_s^\mu)}\grad
\lb(q_s^\mu)\cdot \boldp_s^\mu ds}
\\
\play{=-\int_0^t
e^{-\frac{1}{\mu}A(\mu,s)}\dfrac{1}{\lb^2(q_s^\mu)}e^{-\frac{1}{\mu}A(\mu,s)}\grad
\lb(q_s^\mu) \cdot \left(\boldp+ \dfrac{1}{\mu}\int_0^s
e^{\frac{1}{\mu}A(\mu,r)}\boldb(\boldq_r^\mu)dr +
\frac{1}{\mu}\int_0^s e^{\frac{1}{\mu}A(\mu,r)}d\boldW_r\right) ds}
\\
\play{=(I_1)+(I_2)+(I_3) \ .}
\end{array}$$

Here

$$(I_1)=-\boldp\cdot \int_0^t
e^{-\frac{2}{\mu}A(\mu,s)}\dfrac{\grad
\lb(\boldq_s^\mu)}{\lb^2(\boldq_s^\mu)} ds \ ,$$

$$(I_2)=-\dfrac{1}{\mu}\int_0^t
e^{-\frac{2}{\mu}A(\mu,s)}\dfrac{1}{\lb^2(\boldq_s^\mu)}\grad
\lb(\boldq_s^\mu)\cdot \left(\int_0^s
e^{\frac{1}{\mu}A(\mu,r)}\boldb(\boldq_r^\mu)dr\right) ds \ ,$$

$$(I_3)=-\dfrac{1}{\mu}\int_0^t
e^{-\frac{2}{\mu}A(\mu,s)}\dfrac{1}{\lb^2(q_s^\mu)}\grad
\lb(q_s^\mu)\cdot \left( \int_0^s
e^{\frac{1}{\mu}A(\mu,r)}d\boldW_r\right) ds \ .$$

We can derive, using (2.6) and (2.7), that

$$
|(I_1)|\leq \dfrac{\|\grad
\lb\|_{\infty}}{\lb_0^2}|\boldp|_{\R^d}\int_0^t
e^{-\frac{2}{\mu}\lb_0 s}ds \leq \dfrac{\|\grad
\lb\|_{\infty}}{\lb_0^2}|\boldp|_{\R^d}\dfrac{\mu}{2\lb_0} \ .
$$

$$
\begin{array}{l}
\play{|(I_2)|\leq \dfrac{\|\grad
\lb\|_{\infty}}{\lb_0^2}\|\boldb\|_\infty\dfrac{1}{\mu}\int_0^t
e^{-\frac{2}{\mu}A(\mu,s)}\left(\int_0^s
e^{\frac{1}{\mu}A(\mu,r)}dr\right)ds}
\\
\play{\leq \dfrac{\|\grad
\lb\|_\infty}{\lb_0^2}\|\boldb\|_\infty\dfrac{1}{\mu}\int_0^t
\left(\int_0^s
e^{-\frac{1}{\mu}(s-r)\lb_0}dr\right)e^{-\frac{1}{\mu}\lb_0 s}ds}
\\
\play{= \dfrac{\|\grad
\lb\|_\infty}{\lb_0^2}\|\boldb\|_\infty\dfrac{1}{\mu}\int_0^t
\dfrac{\mu}{\lb_0}(1-e^{-\frac{\lb_0 s}{\mu}})e^{-\frac{\lb_0
s}{\mu}}ds}
\\
\play{\leq \dfrac{\|\grad
\lb\|_\infty}{\lb_0^3}\|\boldb\|_\infty\int_0^t e^{-\frac{\lb_0
s}{\mu}}ds}
\\
\play{\leq \dfrac{\|\grad
\lb\|_\infty}{\lb_0^3}\|\boldb\|_\infty\dfrac{\mu}{\lb_0}} \ .
\end{array}$$

Since

$$
\begin{array}{l}
\play{|(I_3)|\leq \dfrac{\|\grad
\lb\|_{\infty}}{\lb_0^2}\dfrac{1}{\mu}\left|\int_0^t
e^{-\frac{1}{2\mu}A(\mu,s)}\left(\int_0^s
e^{-\frac{1}{2\mu}A(\mu,s)}e^{-\frac{1}{\mu}A(\mu,s)+\frac{1}{\mu}A(\mu,r)}d\boldW_r\right)ds\right|_{\R^d}}
\ ,
\end{array}$$
we could estimate, by Cauchy-Schwarz inequality and (2.6), (2.7),
that

$$
\begin{array}{l}
\play{\E|(I_3)|^2\leq \left(\dfrac{\|\grad
\lb\|_{\infty}}{\lb_0^2}\right)^2\dfrac{1}{\mu^2}\E\left|\int_0^t
e^{-\frac{1}{2\mu}A(\mu,s)}\left(\int_0^s
e^{-\frac{1}{2\mu}A(\mu,s)}e^{-\frac{1}{\mu}A(\mu,s)+\frac{1}{\mu}A(\mu,r)}d\boldW_r\right)ds\right|^2_{\R^d}}
\\
\play{\leq \left(\dfrac{\|\grad
\lb\|_{\infty}}{\lb_0^2}\right)^2\dfrac{1}{\mu^2}\E\left(\int_0^t
e^{-\frac{1}{\mu}A(\mu,s)}ds\right)
\left(\int_0^te^{-\frac{1}{\mu}A(\mu,s)}
\left|\int_0^se^{-\frac{1}{\mu}A(\mu,s)+\frac{1}{\mu}A(\mu,r)}d\boldW_r\right|^2_{\R^d}ds\right)}
\\
\play{\leq \left(\dfrac{\|\grad
\lb\|_{\infty}}{\lb_0^2}\right)^2\dfrac{1}{\mu^2}\left(\int_0^t
e^{-\frac{\lb_0 s}{\mu}}ds\right) \left(\int_0^te^{-\frac{\lb_0
s}{\mu}}\E
\left|\int_0^se^{-\frac{1}{\mu}A(\mu,s)+\frac{1}{\mu}A(\mu,r)}d\boldW_r\right|^2_{\R^d}ds\right)}
\\
\play{= \left(\dfrac{\|\grad
\lb\|_{\infty}}{\lb_0^2}\right)^2\dfrac{1}{\mu^2}\left(\int_0^t
e^{-\frac{\lb_0 s}{\mu}}ds\right) \left(\int_0^te^{-\frac{\lb_0
s}{\mu}} \left(\int_0^s \E
e^{-\frac{2}{\mu}A(\mu,s)+\frac{2}{\mu}A(\mu,r)}dr\right) ds\right)}
\\
\play{\leq \left(\dfrac{\|\grad
\lb\|_{\infty}}{\lb_0^2}\right)^2\dfrac{1}{\mu^2}\left(\int_0^t
e^{-\frac{\lb_0 s}{\mu}}ds\right) \left(\int_0^te^{-\frac{\lb_0
s}{\mu}} \left(\int_0^s e^{-\frac{2\lb_0 s}{\mu}+\frac{2\lb_0
r}{\mu}}dr\right) ds\right)}
\\
\play{\leq \left(\dfrac{\|\grad
\lb\|_{\infty}}{\lb_0^2}\right)^2\dfrac{1}{\mu^2}\dfrac{\mu}{\lb_0}
\left(\int_0^te^{-\frac{\lb_0 s}{\mu}} \dfrac{\mu}{2\lb_0}
ds\right)}
\\
\play{\leq \left(\dfrac{\|\grad
\lb\|_{\infty}}{\lb_0^2}\right)^2\dfrac{\mu}{2\lb_0^3}} \ .
\end{array}$$

Combining these estimates we see that $\E|(I)|^2\ra 0$ as $\mu \da
0$.

So $\E|\boldal(\mu)|_{\R^d}^2\ra 0$ as $\mu\da 0$ for any
$|\boldp|_{\R^d}<\infty$.

\

\textbf{2.2. Estimates of $\boldbt(\mu)$. }

\

We have, by (2.5), that

$$\begin{array}{l}
\play{\boldbt(\mu)=\dfrac{1}{\mu}\int_0^t
e^{-\frac{1}{\mu}A(\mu,s)}\left(\int_0^s
e^{\frac{1}{\mu}A(\mu,r)}\boldb(\boldq_r^\mu)dr\right)ds}
\\
\play{=\dfrac{1}{\mu}\int_0^t \left(\int_0^s
e^{\frac{1}{\mu}A(\mu,r)}\boldb(\boldq_r^\mu)dr\right)\left(-\dfrac{\mu}{\lb(\boldq_s^\mu)}\right)d(e^{-\frac{1}{\mu}A(\mu,s)})}
\\
\play{=\int_0^t \left(\int_0^s
e^{\frac{1}{\mu}A(\mu,r)}\boldb(\boldq_r^\mu)dr\right)\left(-\dfrac{1}{\lb(\boldq_s^\mu)}\right)d(e^{-\frac{1}{\mu}A(\mu,s)})}
\\
\play{=-\dfrac{e^{-\frac{1}{\mu}A(\mu,s)}}{\lb(\boldq_s^\mu)}\left.\int_0^s
e^{\frac{1}{\mu}A(\mu,r)}\boldb(\boldq_r^\mu)dr\right|^{s=t}_{s=0}+\int_0^t
e^{-\frac{1}{\mu}A(\mu,s)}d\left(\int_0^s
e^{\frac{1}{\mu}A(\mu,r)}\boldb(\boldq_r^\mu)dr\dfrac{1}{\lb(\boldq_s^\mu)}\right)}
\\
\play{=-\dfrac{e^{-\frac{1}{\mu}A(\mu,t)}}{\lb(\boldq_t^\mu)}\int_0^t
e^{\frac{1}{\mu}A(\mu,s)}\boldb(\boldq_s^\mu)ds+\int_0^t
\dfrac{\boldb(\boldq_s^\mu)}{\lb(\boldq_s^\mu)}ds+\int_0^t
e^{-\frac{1}{\mu}A(\mu,s)}\left(\int_0^s
e^{\frac{1}{\mu}A(\mu,r)}\boldb(\boldq_r^\mu)dr\right)d\left(\frac{1}{\lb(\boldq_s^\mu)}\right)}
\\
\play{=\boldR_{\boldbt}(\mu)+\int_0^t
\dfrac{\boldb(\boldq_s^\mu)}{\lb(\boldq_s^\mu)}ds+\vec{(II)} \ .}
\end{array}$$

It is easy to see that

$$\begin{array}{l}
\play{|\boldR_{\boldbt}(\mu)|_{\R^d}\leq\dfrac{\|\boldb\|_\infty}{\lb_0}\int_0^t
e^{-\frac{\lb_0}{\mu}(t-s)}ds=\dfrac{\|\boldb\|_\infty}{\lb_0}\dfrac{\mu}{\lb_0}}
\ .
\end{array}$$

We also have

$$\begin{array}{l}
\play{\vec{(II)}=-\int_0^t e^{-\frac{1}{\mu}A(\mu,s)}\left(\int_0^s
e^{\frac{1}{\mu}A(\mu,r)}\boldb(\boldq_r^\mu)dr\right)\dfrac{1}{\lb^2(\boldq_s^\mu)}\grad
\lb(\boldq_s^\mu)\cdot \boldp_s^\mu ds}
\\
\play{=-\int_0^t e^{-\frac{1}{\mu}A(\mu,s)}\left(\int_0^s
e^{\frac{1}{\mu}A(\mu,r)}\boldb(\boldq_r^\mu)dr\right)\dfrac{1}{\lb^2(\boldq_s^\mu)}\grad
\lb(\boldq_s^\mu)\cdot}
\\
\play{ \ \ \ \ \ \ \ \ \ \ \cdot
e^{-\frac{1}{\mu}A(\mu,s)}\left(\boldp+ \dfrac{1}{\mu}\int_0^s
e^{\frac{1}{\mu}A(\mu,r)}\boldb(\boldq_r^\mu)dr +
\frac{1}{\mu}\int_0^s e^{\frac{1}{\mu}A(\mu,r)}d\boldW_r\right)ds}
\\
\play{=\vec{(II_1)}+\vec{(II_2)}+\vec{(II_3)} \ .}
\end{array}$$

Here

$$\vec{(II_1)}=-\int_0^t \dfrac{e^{-\frac{2}{\mu}A(\mu,s)}}{\lb^2(\boldq_s^\mu)}
\left(\int_0^s
e^{\frac{1}{\mu}A(\mu,r)}\boldb(\boldq_r^\mu)dr\right)\grad
\lb(\boldq_s^\mu)\cdot \boldp ds \ , $$

$$\vec{(II_2)}=-\dfrac{1}{\mu}\int_0^t \dfrac{e^{-\frac{2}{\mu}A(\mu,s)}}{\lb^2(\boldq_s^\mu)}
\left(\int_0^s
e^{\frac{1}{\mu}A(\mu,r)}\boldb(\boldq_r^\mu)dr\right)\grad
\lb(\boldq_s^\mu)\cdot \left(\int_0^s
e^{\frac{1}{\mu}A(\mu,r)}\boldb(\boldq_r^\mu)dr\right) ds \ ,$$

$$\vec{(II_3)}=-\dfrac{1}{\mu}\int_0^t \dfrac{e^{-\frac{2}{\mu}A(\mu,s)}}{\lb^2(\boldq_s^\mu)}
\left(\int_0^s
e^{\frac{1}{\mu}A(\mu,r)}\boldb(\boldq_r^\mu)dr\right)\grad
\lb(\boldq_s^\mu)\cdot \left(\int_0^s
e^{\frac{1}{\mu}A(\mu,r)}d\boldW_r\right) ds \ .$$

We conclude that

$$\begin{array}{l}
\play{|\vec{(II_1)}|_{\R^d}\leq\dfrac{\|\grad
\lb\|_\infty}{\lb_0^2}|\boldp|_{\R^d}\|\boldb\|_\infty\int_0^t
e^{-\frac{\lb_0 s}{\mu}}\left(\int_0^s
e^{-\frac{(s-r)\lb_0}{\mu}}dr\right)ds}
\\
\play{\leq\dfrac{\|\grad
\lb\|_\infty}{\lb_0^2}|\boldp|_{\R^d}\|\boldb\|_\infty\dfrac{\mu^2}{\lb_0^2}}
\ ;
\end{array}$$

$$\begin{array}{l}
\play{|\vec{(II_2)}|_{\R^d}\leq\dfrac{1}{\mu}\dfrac{\|\grad
\lb\|_\infty}{\lb_0^2}\|\boldb\|_\infty^2\int_0^t \left(\int_0^s
e^{-\frac{(s-r)\lb_0}{\mu}}dr\right)^2ds}
\\
\play{\leq\dfrac{\|\grad
\lb\|_\infty}{\lb_0^2}\|\boldb\|_\infty^2\dfrac{\mu t}{\lb_0^2}} \ ;
\end{array}$$

$$\begin{array}{l}
\play{\E|\vec{(II_3)}|^2_{\R^d}\leq\left(\dfrac{1}{\mu}\dfrac{\|\grad
\lb\|_\infty}{\lb_0^2}\|\boldb\|_\infty\right)^2\E\left|\int_0^t
e^{-\frac{2}{\mu}A(\mu,s)}\left(\int_0^se^{\frac{1}{\mu}A(\mu,r)}dr\right)
\left(\int_0^s
e^{\frac{1}{\mu}A(\mu,r)}d\boldW_r\right)ds\right|^2_{\R^d}}
\\
\play{=\left(\dfrac{1}{\mu}\dfrac{\|\grad
\lb\|_\infty}{\lb_0^2}\|\boldb\|_\infty\right)^2\E\left|\int_0^t
\left(\int_0^se^{-\frac{1}{\mu}(A(\mu,s)-A(\mu,r))}dr\right)
\left(\int_0^s
e^{-\frac{1}{\mu}(A(\mu,s)-A(\mu,r))}d\boldW_r\right)ds\right|^2_{\R^d}}
\\
\play{\leq\left(\dfrac{1}{\mu}\dfrac{\|\grad
\lb\|_\infty}{\lb_0^2}\|\boldb\|_\infty\right)^2\E\left(\int_0^t
\left(\int_0^se^{-\frac{1}{\mu}(A(\mu,s)-A(\mu,r))}dr\right)^2ds\right)
\left(\int_0^t\left|\int_0^s
e^{-\frac{1}{\mu}(A(\mu,s)-A(\mu,r))}d\boldW_r\right|^2_{\R^d}ds\right)}
\\
\play{\leq\left(\dfrac{1}{\mu}\dfrac{\|\grad
\lb\|_\infty}{\lb_0^2}\|\boldb\|_\infty\right)^2\left(\int_0^t
\left(\int_0^se^{-\frac{(s-r)\lb_0}{\mu}}dr\right)^2ds\right)
\left(\int_0^t\E\left|\int_0^s
e^{-\frac{1}{\mu}(A(\mu,s)-A(\mu,r))}d\boldW_r\right|^2_{\R^d}ds\right)}
\\
\play{\leq\left(\dfrac{1}{\mu}\dfrac{\|\grad
\lb\|_\infty}{\lb_0^2}\|\boldb\|_\infty\right)^2\left(\int_0^t
\left(\int_0^se^{-\frac{(s-r)\lb_0}{\mu}}dr\right)^2ds\right)
\left(\int_0^t\left(\int_0^s
e^{-\frac{2(s-r)\lb_0}{\mu}}dr\right)ds\right)}
\\
\play{\leq\left(\dfrac{\|\grad
\lb\|_\infty}{\lb_0^2}\|\boldb\|_\infty\right)^2\left(\dfrac{
t}{\lb_0}\right)^2 \left(\dfrac{\mu t}{2\lb_0}\right) \ .}
\end{array}$$

Combining these estimates we see that $\E|\vec{(II)}|^2_{\R^d}\ra 0$
as $\mu \da 0$. This implies that
$\play{\E\left|\boldbt(\mu)-\int_0^t\dfrac{\boldb(\boldq_s^\mu)}{\lb(\boldq_s^\mu)}ds\right|_{\R^d}^2\ra
0}$ as $\mu \da 0$.

\

\textbf{2.3. Estimates of $\boldgm(\mu)$ - the reason why the
classical Smoluchowski-Kramers approximation does not work.}

\

We will show that $\E\left|\boldgm(\mu)-\play{\int_0^t
\dfrac{1}{\lb(\boldq_s^\mu)}d\boldW_s}\right|_{\R^d}^2$, in general,
does not tend to $0$ as $\mu \da 0$. Therefore the
Smoluchowski-Kramers approximation does not work in the case of
purely white noise perturbation.

We have, by (2.5), that

$$\begin{array}{l}
\play{\boldgm(\mu)=\dfrac{1}{\mu}\int_0^t
e^{-\frac{1}{\mu}A(\mu,s)}\left(\int_0^s
e^{\frac{1}{\mu}A(\mu,r)}d\boldW_r\right)ds}
\\
\play{=\dfrac{1}{\mu}\int_0^t \left(\int_0^s
e^{\frac{1}{\mu}A(\mu,r)}d\boldW_r\right)\left(
-\dfrac{\mu}{\lb(\boldq_s^\mu)}\right)d(e^{-\frac{1}{\mu}A(\mu,s)})}
\\
\play{=-\left[\dfrac{\play{\int_0^t
e^{\frac{1}{\mu}A(\mu,s)}d\boldW_s}}{\lb(\boldq_t^\mu)}e^{-\frac{1}{\mu}A(\mu,t)}
-\int_0^t
e^{-\frac{1}{\mu}A(\mu,s)}d\left(\dfrac{1}{\lb(\boldq_s^\mu)}\int_0^s
e^{\frac{1}{\mu}A(\mu,r)}d\boldW_r\right)\right]}
\\
\play{=-\dfrac{\play{\int_0^t
e^{\frac{1}{\mu}A(\mu,s)}d\boldW_s}}{\lb(\boldq_t^\mu)}e^{-\frac{1}{\mu}A(\mu,t)}
+ \int_0^t \dfrac{1}{\lb(\boldq_s^\mu)}d\boldW_s + }
\\
\play{ \ \ \ \ \ \ \ \ \ \ \ \ \ + \int_0^t
e^{-\frac{1}{\mu}A(\mu,s)}\left(\int_0^s
e^{\frac{1}{\mu}A(\mu,r)}d\boldW_r\right)d\left(\dfrac{1}{\lb(\boldq_s^\mu)}\right)}
\\
\play{=\boldR_{\boldgm}(\mu)+\int_0^t
\dfrac{1}{\lb(\boldq_s^\mu)}d\boldW_s+\vec{(III)} \ .}
\end{array}$$

It is easy to check that

$$\E|\boldR_{\boldgm}(\mu)|^2_{\R^d}\leq \dfrac{1}{\lb_0^2}
\int_0^t e^{-\frac{2\lb_0(t-s)}{\mu}}ds\leq \dfrac{\mu}{2\lb_0^3} \
.$$

We have

$$\begin{array}{l}
\play{\vec{(III)}=\int_0^t e^{-\frac{1}{\mu}A(\mu,s)}\left(\int_0^s
e^{\frac{1}{\mu}A(\mu,r)}d\boldW_r\right)\left(-\dfrac{1}{\lb^2(\boldq_s^\mu)}\right)\grad
\lb(\boldq_s^\mu)\cdot \boldp_s^\mu ds}
\\
\play{=\vec{(III_1)}+\vec{(III_2)}+\vec{(III_3)}} \end{array}$$
where

$$\vec{(III_1)}=-\int_0^t \dfrac{e^{-\frac{2}{\mu}A(\mu,s)}}{\lb^2(\boldq_s^\mu)}
\left(\int_0^s e^{\frac{1}{\mu}A(\mu,r)}d\boldW_r\right)\grad
\lb(\boldq_s^\mu)\cdot \boldp ds \ ,$$

$$\vec{(III_2)}=-\dfrac{1}{\mu}\int_0^t \dfrac{e^{-\frac{2}{\mu}A(\mu,s)}}{\lb^2(\boldq_s^\mu)}
\left(\int_0^s e^{\frac{1}{\mu}A(\mu,r)}d\boldW_r\right)\grad
\lb(\boldq_s^\mu)\cdot \left(\int_0^s
e^{\frac{1}{\mu}A(\mu,r)}\boldb(\boldq_r^\mu)dr\right) ds \ ,$$

$$\vec{(III_3)}=-\dfrac{1}{\mu}\int_0^t \dfrac{e^{-\frac{2}{\mu}A(\mu,s)}}{\lb^2(\boldq_s^\mu)}
\left(\int_0^s e^{\frac{1}{\mu}A(\mu,r)}d\boldW_r\right)\grad
\lb(\boldq_s^\mu)\cdot \left(\int_0^s
e^{\frac{1}{\mu}A(\mu,r)}d\boldW_r\right) ds \ ,$$

We can estimate

$$\begin{array}{l}
\play{\E|\vec{(III_1)}|^2_{\R^d}\leq
\left(\dfrac{|\boldp|_{\R^d}\|\grad \lb\|_\infty}{\lb_0^2}\right)^2
\E\left|\int_0^t
e^{-\frac{1}{\mu}A(\mu,s)}\left(\int_0^se^{-\frac{1}{\mu}(A(\mu,s)-A(\mu,r))}d\boldW_r\right)ds\right|^2_{\R^d}}
\\
\play{\leq \left(\dfrac{|\boldp|_{\R^d}\|\grad
\lb\|_\infty}{\lb_0^2}\right)^2 \E\left(\int_0^t
e^{-\frac{2}{\mu}A(\mu,s)}ds\right)
\left(\int_0^t\left|\int_0^se^{-\frac{1}{\mu}(A(\mu,s)-A(\mu,r))}d\boldW_r\right|^2_{\R^d}ds\right)}
\\
\play{\leq \left(\dfrac{|\boldp|_{\R^d}\|\grad
\lb\|_\infty}{\lb_0^2}\right)^2 \left(\int_0^t e^{-\frac{2\lb_0
s}{\mu}}ds\right) \left(\int_0^t\left(\int_0^se^{-\frac{2\lb_0
(s-r)}{\mu}}dr\right)ds\right)}
\\
\play{\leq \left(\dfrac{|\boldp|_{\R^d}\|\grad
\lb\|_\infty}{\lb_0^2}\right)^2 \left(\dfrac{\mu}{2\lb_0}\right)
\left(\dfrac{\mu t}{2\lb_0}\right)} \ .
\end{array}$$

The term $\vec{(III_2)}$ could be estimated in the same way as
$\vec{(II_3)}$:

$$\begin{array}{l}
\play{\E|\vec{(III_2)}|^2_{\R^d}\leq\left(\dfrac{\|\grad
\lb\|_\infty}{\lb_0^2}\|\boldb\|_\infty\right)^2\left(\dfrac{
t}{\lb_0}\right)^2 \left(\dfrac{\mu t}{2\lb_0}\right) \ .}
\end{array}$$

But in general one cannot estimate $\E|\vec{(III_3)}|^2$ up to a
term which goes to $0$ as $\mu \da 0$. As an example, let
$\Lb=\|\lb\|_\infty$ and let us suppose that for $0\leq t \leq
T<\infty$ we have $\grad \lb(\boldq_t^\mu)=\mathbf{e}_1$. Here
$\mathbf{e}_1$ is the unit basis vector $\mathbf{e}_1=(1,0,...,0)$
in $\R^d$. Let $W_r^k$ be the $k$-th ($1\leq k \leq d$) component of
the Wiener process $\boldW_r$. We have, for $0<t\leq T$:

$$\begin{array}{l}
\play{\E|\vec{(III_3)}|_{\R^d}\geq \dfrac{1}{\mu
\Lb^2}\E\left|\int_0^t\left( \int_0^s
e^{-\frac{1}{\mu}(A(\mu,s)-A(\mu,r))}d\boldW_r\right)\left( \int_0^s
e^{-\frac{1}{\mu}(A(\mu,s)-A(\mu,r))}dW_r^1\right)ds\right|_{\R^d}}
\\
\play{=\dfrac{1}{\mu \Lb^2}\E \left[\left(\int_0^t\left(\int_0^s
e^{-\frac{1}{\mu}(A(\mu,s)-A(\mu,r))}dW_r^1\right)^2ds\right)^2
+\right.}
\\
\play{\ \ \ \ \ \ \ \ \ \ \
 \ \left.+\sum\li_{k=2}^d\left(\int_0^t\left(\int_0^se^{-\frac{1}{\mu}(A(\mu,s)-A(\mu,r))}dW_r^k\right)
 \left(\int_0^se^{-\frac{1}{\mu}(A(\mu,s)-A(\mu,r))}dW_r^1\right)ds\right)^2\right]^{\frac{1}{2}}}
\\
\play{\geq \dfrac{1}{\mu \Lb^2}\E\left(\int_0^t\left( \int_0^s
e^{-\frac{1}{\mu}(A(\mu,s)-A(\mu,r))}dW_r^1\right)^2ds\right)}
\\
\play{= \dfrac{1}{\mu \Lb^2}\left(\int_0^t \left(\int_0^s \E
e^{-\frac{2}{\mu}(A(\mu,s)-A(\mu,r))}dr\right)ds\right)}
\\
\play{\geq \dfrac{1}{\mu \Lb^2}\left(\int_0^t \left(\int_0^s
e^{-\frac{2}{\mu}\Lb(s-r)}dr\right)ds\right)}
\\
\play{=\dfrac{1}{\mu\Lb^2}\dfrac{\mu}{2\Lb}\int_0^t(1-e^{-\frac{2\Lb
s}{\mu}})ds=\dfrac{t}{2\Lb^3}-\dfrac{\mu}{4\Lb^4}(1-e^{-\frac{2\Lb
t}{\mu}})} \ ,
\end{array}$$
which does not tend to $0$ as $\mu \da 0$. Since
$\E|\vec{(III_3)}|_{\R^d}^2\geq (\E|\vec{(III_3)}|_{\R^d})^2$, we
see that $\E|\vec{(III_3)}|_{\R^d}^2$ does not go to $0$ as $\mu \da
0$. Now we have

$$\E\left|\boldgm(\mu)-\play{\int_0^t
\dfrac{1}{\lb(\boldq_s^\mu)}d\boldW_s}\right|_{\R^d}^2\geq
\dfrac{1}{4}\E|\vec{(III_3)}|_{\R^d}^2-\E|\boldR_{\boldgm}(\mu)|_{\R^d}^2-
\E|\vec{(III_1)}|_{\R^d}^2-\E|\vec{(III_2)}|_{\R^d}^2 \ .$$
Therefore $\play{\E\left|\boldgm(\mu)-\play{\int_0^t
\dfrac{1}{\lb(\boldq_s^\mu)}d\boldW_s}\right|_{\R^d}^2}$ is
uniformly bounded from below by a positive constant as $\mu \da 0$.

\

We can check now that the process $\boldq_t^\mu$, $0\leq t \leq T$,
does not converge as $\mu \da 0$ to the process $\boldq_t$,
$\boldq_0=\boldq$. We have

$$\begin{array}{l}
\play{\boldq_t^\mu=\boldq+\int_0^t\dfrac{\boldb(\boldq_s^\mu)}{\lb(\boldq_s^\mu)}ds+
\int_0^t\dfrac{1}{\lb(\boldq_s^\mu)}d\boldW_s+}
\\
\play{\ \ \ \ \ \ \ \ \ \ \ \ \ \ \ \ \
+\boldal(\mu)+\left(\boldbt(\mu)-\int_0^t\dfrac{\boldb(\boldq_s^\mu)}{\lb(\boldq_s^\mu)}ds\right)
+\left(\boldgm(\mu)-\int_0^t\dfrac{1}{\lb(\boldq_s^\mu)}d\boldW_s\right)}
\ ,
\\
\play{\boldq_t=\boldq+\int_0^t\dfrac{\boldb(\boldq_s)}{\lb(\boldq_s)}ds+
\int_0^t\dfrac{1}{\lb(\boldq_s)}d\boldW_s} \ .
\end{array}$$

Suppose that we have, for any $\kp, T>0$ and any
$\boldp_0^\mu=\boldp\in \R^d$ fixed, that $$\lim\li_{\mu \da
0}\Prob\left(\max\li_{0\leq t \leq T}|\boldq_t^\mu-
\boldq_t|_{\R^d}^2\geq \kp\right)=0 \ .$$

We have, for some $A>0$ independent of $\mu$ and $\kp$, that

$$\begin{array}{l}
\play{\E\left|(\boldq_t^\mu-\boldq_t)-\int_0^t\left(\dfrac{\boldb(\boldq_s^\mu)}{\lb(\boldq_s^\mu)}
-\dfrac{\boldb(\boldq_s)}{\lb(\boldq_s)}\right)ds-\int_0^t
\left(\dfrac{1}{\lb(\boldq_s^\mu)}-\dfrac{1}{\lb(\boldq_s)}\right)d\boldW_s\right|_{\R^d}^2}
\\
\play{\leq A\E\max\li_{0\leq s\leq
t}|\boldq_s^\mu-\boldq_s|_{\R^d}^2}
\\
\play{\leq A \left[\Prob\left(\max\li_{0\leq s \leq t}|\boldq_s^\mu-
\boldq_s|_{\R^d}^2\geq \kp\right)\cdot\E\max\li_{0\leq s\leq
t}|\boldq_s^\mu-\boldq_s|_{\R^d}^2+\Prob\left(\max\li_{0\leq s \leq
t}|\boldq_s^\mu- \boldq_s|_{\R^d}^2< \kp\right)\cdot\kp\right]}
\\
\play{\leq A[\kp+o(\mu,\kp)] \ ,}
\end{array}$$
since $\E\max\li_{0\leq s\leq
t}|\boldq_s^\mu-\boldq_s|_{\R^d}^2<\infty$. Here the term
$o(\mu,\kp)$ converges to $0$ as $\mu \da 0$ for every fixed $\kp
>0$. Fix $\kp
>0$, let $\mu \da 0$, we see that

$$\lim\li_{\mu \da 0}\E\left|(\boldq_t^\mu-\boldq_t)-\int_0^t\left(\dfrac{\boldb(\boldq_s^\mu)}{\lb(\boldq_s^\mu)}
-\dfrac{\boldb(\boldq_s)}{\lb(\boldq_s)}\right)ds-\int_0^t
\left(\dfrac{1}{\lb(\boldq_s^\mu)}-\dfrac{1}{\lb(\boldq_s)}\right)d\boldW_s\right|_{\R^d}^2\leq
A \kp \ .$$

Since $\kp>0$ is arbitrary, we see that

$$\lim\li_{\mu \da 0}\E\left|(\boldq_t^\mu-\boldq_t)-\int_0^t\left(\dfrac{\boldb(\boldq_s^\mu)}{\lb(\boldq_s^\mu)}
-\dfrac{\boldb(\boldq_s)}{\lb(\boldq_s)}\right)ds-\int_0^t
\left(\dfrac{1}{\lb(\boldq_s^\mu)}-\dfrac{1}{\lb(\boldq_s)}\right)d\boldW_s\right|_{\R^d}^2=0
\ .$$

On the other hand, let us suppose that $\grad
\lb(\boldq_t^\mu)=\mathbf{e}_1$ for $0\leq t \leq T<\infty$. Here
$\mathbf{e}_1$ is the unit basis vector $\mathbf{e}_1=(1,0,...,0)$
in $\R^d$. We have

$$\begin{array}{l}
\play{\E\left|\boldal(\mu)+\left(\boldbt(\mu)-\int_0^t\dfrac{\boldb(\boldq_s^\mu)}{\lb(\boldq_s^\mu)}ds\right)
+\left(\boldgm(\mu)-\int_0^t\dfrac{1}{\lb(\boldq_s^\mu)}d\boldW_s\right)\right|_{\R^d}^2}
\\
\play{\geq
\dfrac{1}{3}\E\left|\boldgm(\mu)-\int_0^t\dfrac{1}{\lb(\boldq_s^\mu)}d\boldW_s\right|_{\R^d}^2
-\E\left|\boldal(\mu)\right|_{\R^d}^2
-\E\left|\boldbt(\mu)-\int_0^t\dfrac{\boldb(\boldq_s^\mu)}{\lb(\boldq_s^\mu)}ds\right|_{\R^d}^2
\ .}
\end{array}$$

It follows from our estimates that this leads to a contradiction.

\

\section{Regularization via approximation of the Wiener process}

We could regularize the problem via \textit{approximation of the
Wiener process}. To this end we introduce the process

$$\boldW_t^\dt=\dfrac{1}{\dt}\int_0^\infty \boldW_s \rho\left(
\dfrac{s-t}{\dt}\right)ds=\dfrac{1}{\dt}\int_0^\dt
\boldW_{s+t}\rho\left(\dfrac{s}{\dt}\right)ds \ ,$$ where
$\rho(\bullet)$ is a smooth $C^\infty$ function whose support is
contained in the interval $[0,1]$ such that

$$\int_0^1 \rho(s) ds =1 \ .$$

One can prove that (see [3] and the references there)

$$\lim\li_{\dt\da 0}\E\max\li_{t\in [0,T]}|\boldW_t^\dt-\boldW_t|^2_{\R^d}=0 \ .$$

We have

$$\dot{\boldW}_t^\dt=-\dfrac{1}{\dt}\int_0^1 \boldW_{t+\dt r}\dot{\rho}(r)dr \ .$$

We can then introduce the following regularization of our problem:
first we consider the system

$$\mu \ddot{\boldq}^{\mu,\dt}_t=\boldb(\boldq_t^{\mu,\dt})-\lb(\boldq_t^{\mu,\dt})\dot{\boldq}_t^{\mu,\dt}
+\dot{\boldW}_t^\dt \ , \ \boldq_0^{\mu,\dt}=\boldq\in \R^d
 \ , \ \dot{\boldq}_0^{\mu,\dt}=\boldp\in \R^d. \eqno(3.1)$$

Equivalently it is the first order system

$$
\left\{\begin{array}{l} \dot{\boldq}_t^{\mu,\dt}=\boldp_t^{\mu,\dt}  \ , \\
\dot{\boldp}_t^{\mu,\dt}=\dfrac{1}{\mu}\boldb(\boldq_t^{\mu,\dt})-\dfrac{\lb(\boldq_t^{\mu,\dt})}{\mu}\boldp_t^{\mu,\dt}
+\dfrac{1}{\mu}\dot{\boldW}_t^\dt \ .
\end{array}
\right. \eqno(3.2)$$

We can proceed with the estimates similar to in the previous
sections. Since for fixed $\dt>0$ ,
$$|\dot{\boldW}_t^\dt|_{\R^d}\leq \dfrac{1}{\dt}
\left(\max\li_{0\leq r \leq 1}|\dot{\rho}(r)|\right)
\left(\max\li_{t\leq s \leq t+\dt}|\boldW_s|_{\R^d}\right)<\infty \
, \text{ a. s. } \ , \eqno(3.3)$$ we could prove that all the terms
$$\E|\boldal(\mu)|_{\R^d}
\ , \
\E\left|\play{\boldbt(\mu)-\int_0^t\dfrac{\boldb(\boldq_s^\mu)}{\lb(\boldq_s^\mu)}ds}\right|_{\R^d}
\ , \
\E\left|\play{\boldgm(\mu)-\int_0^t\dfrac{1}{\lb(\boldq_s^\mu)}d\boldW_s^\dt}\right|_{\R^d}$$
goes to zero as $\mu \da 0$. (To be precise, we should write
$\boldal(\mu,\dt)$, $\boldbt(\mu,\dt)$ and $\boldgm(\mu,\dt)$ to
indicate the dependence on $\dt$, but for brevity we neglect that.)
In particular, with $\dt>0$ fixed, we can estimate the term
$\vec{(III_3)}$ up to a term which tends to $0$ as $\mu \da 0$. We
have

$$\begin{array}{l}
\play{\E|\vec{(III_3)}|_{\R^d}\leq \dfrac{1}{\mu} \dfrac{\|\grad
\lb\|_\infty}{\lb_0^2}\int_0^t \E\left|\int_0^s
e^{-\frac{1}{\mu}(A(\mu,s)-A(\mu,r))}\dot{\boldW}_r^\dt
dr\right|^2_{\R^d}ds}
\\
\play{= \dfrac{1}{\mu} \dfrac{\|\grad \lb\|_\infty}{\lb_0^2}\int_0^t
\dfrac{1}{\dt^2}\E\left|\int_0^s
e^{-\frac{1}{\mu}(A(\mu,s)-A(\mu,r))}\left(\int_0^1 \boldW_{r+\dt
m}\dot{\rho}(m)dm\right)dr\right|^2_{\R^d}ds}
\\
\play{= \dfrac{1}{\mu} \dfrac{\|\grad \lb\|_\infty}{\lb_0^2}\int_0^t
\dfrac{1}{\dt^2}\E\left|\int_0^1 \dot{\rho}(m)\boldW_{r+\dt m}dm
\int_0^s e^{-\frac{1}{\mu}(A(\mu,s)-A(\mu,r))}dr\right|^2_{\R^d}ds}
\\
\play{\leq \dfrac{1}{\mu} \dfrac{\|\grad
\lb\|_\infty}{\lb_0^2}\int_0^t \dfrac{1}{\dt^2}\left(\max\li_{0\leq
m \leq 1}|\dot{\rho}(m)|\right)^2\E\left(\max\li_{0\leq l \leq
s+\dt}|\boldW_l|_{\R^d}\right)^2\left(\int_0^s
e^{-\frac{\lb_0(s-r)}{\mu}}dr\right)^2ds}
\\
\play{\leq \mu \dfrac{\|\grad \lb\|_\infty}{\lb_0^4}
\dfrac{t}{\dt^2}\left(\max\li_{0\leq m \leq
1}|\dot{\rho}(m)|\right)^2\E\left(\max\li_{0\leq l \leq
s+\dt}|\boldW_l|_{\R^d}\right)^2} \ .
\end{array}$$

Therefore, for fixed $\dt > 0$, we have $\E|\vec{(III_3)}|_{\R^d}\ra
0$ as $\mu \da 0$. By (2.4), we get:

$$\begin{array}{l}
\play{\boldq_t^{\mu,\dt}=\boldq+\int_0^t\dfrac{\boldb(\boldq_s^{\mu,\dt})}{\lb(\boldq_s^{\mu,\dt})}ds+
\int_0^t\dfrac{1}{\lb(\boldq_s^{\mu,\dt})}d\boldW_s^\dt+}
\\
\play{\ \ \ \ \ \ \ \ \ \ \ \ \ \ \ \ \
+\boldal(\mu)+\left(\boldbt(\mu)-\int_0^t\dfrac{\boldb(\boldq_s^{\mu,\dt})}{\lb(\boldq_s^{\mu,\dt})}ds\right)
+\left(\boldgm(\mu)-\int_0^t\dfrac{1}{\lb(\boldq_s^{\mu,\dt})}d\boldW_s^\dt\right)}
\ .
\end{array} \eqno(3.4)$$

Let the process $\widetilde{\boldq}_t^\dt$ be governed by the
equation

$$\dot{\widetilde{\boldq}}_t^\dt=\dfrac{\boldb(\widetilde{\boldq}_t^\dt)}{\lb(\widetilde{\boldq}_t^\dt)}
+\dfrac{1}{\lb(\widetilde{\boldq}_t^\dt)} \dot{\boldW}_t^\dt \ , \
\widetilde{\boldq}_0^\dt=\boldq\in \R^d \ . \eqno(3.5)$$

Then

$$\tilde{\boldq}_t^\dt=\boldq+ \int_0^t \dfrac{\boldb(\tilde{\boldq}_s^\dt)}{\lb(\tilde{\boldq}_s^\dt)}ds+
\int_0^t \dfrac{1}{\lb(\tilde{\boldq}_s^\dt)}d\boldW_s^\dt \  .
\eqno(3.6)$$

Let $M(t,\dt,\mu)=\E\max\li_{0\leq s\leq
t}|\boldq_s^{\mu,\dt}-\tilde{\boldq}_s^\dt|_{\R^d}$. By (3.4) and
(3.6), using estimate (3.3), we have

$$M(t,\dt,\mu)\leq K_1\int_0^t M(s,\dt,\mu)ds + K_2(t,\dt)\int_0^t
M(s,\dt,\mu)ds + o_\mu(1) \ . $$

Here $o_\mu(1)$ is a term which goes to 0 as $\mu \da 0$. The
positive constant $K_1$ is independent of $\mu$, $\dt$ and $t$. The
positive constant $K_2=K_2(t,\dt)$ may depend on $t$ and $\dt$, but
is independent of $\mu$. Now we use the Bellman-Gronwall inequality:

$$M(t,\dt,\mu)\leq o_\mu(1)\exp((K_1+K_2(t,\dt))t) \ .$$

We conclude that for any $\dt ,\kp , T
>0$ fixed and any $\boldp_0^{\mu,\dt}=\boldp$ fixed,

$$\lim\li_{\mu \da 0}\Prob\left(\max\li_{0\leq t \leq T}|\boldq_t^{\mu,\dt}-\widetilde{\boldq}_t^\dt|_{\R^d}>\kp\right)=0 \ .$$

Now we can take $\dt \da 0$. Using Theorem 6.7.2 from [10] we get
the following result.

\

\textbf{Theorem 3.1.}\textit{ We have, as $\dt \da 0$, that
$$\lim\li_{\dt \ra 0}\E \max\li_{t\in
[0,T]}|\widetilde{\boldq}_t^\dt-\widehat{\boldq}_t|_{\R^d}=0  \ , $$
where $\widehat{\boldq}_t$ is the solution of the problem}

$$\dot{\widehat{\boldq}}_t=\dfrac{\boldb(\widehat{\boldq}_t)}{\lb(\widehat{\boldq}_t)}+
\dfrac{1}{\lb(\widehat{\boldq}_t)}\circ \dot{\boldW}_t \ , \
\widehat{\boldq}_0=\boldq\in \R^d \ . \eqno(3.6)$$

\

Here the stochastic term is understood in the Stratonovich sense.

\

In the general case

$$\mu \ddot{\boldq}_t^{\mu,\dt}=\boldb(\boldq_t^{\mu,\dt})-
\lb(\boldq_t^{\mu,\dt})\dot{\boldq}_t^{\mu,\dt}
+\sm(\boldq_t^{\mu,\dt})\dot{\boldW}_t^\dt \ , \
\boldq_0^{\mu,\dt}=\boldq \ , \ \dot{\boldq}_0^{\mu,\dt}=\boldp \ ,
\eqno(3.7)$$ where the matrix $\sm(\bullet)$ satisfy assumptions
made in Section 1, we have, similarly, that for any $\dt,\kp , T
>0$ fixed and any $\boldp_0^{\mu,\dt}=\boldp$ fixed,

$$\lim\li_{\mu \da 0}\Prob\left(\max\li_{0\leq t \leq T}|\boldq_t^{\mu,\dt}-\widetilde{\boldq}_t^\dt|_{\R^d}>\kp\right)=0 \ .$$

The process $\widetilde{\boldq}_t^\dt$ is governed by the equation

$$\dot{\widetilde{\boldq}}_t^\dt=\dfrac{\boldb(\widetilde{\boldq}_t^\dt)}{\lb(\widetilde{\boldq}_t^\dt)}
+\dfrac{\sm(\widetilde{\boldq}_t^\dt)}{\lb(\widetilde{\boldq}_t^\dt)}
\dot{\boldW}_t^\dt \ , \ \widetilde{\boldq}_0^\dt=\boldq\in \R^d \ .
\eqno(3.8)$$

And we conclude with

\

\textbf{Theorem 3.2.}\textit{ Under the assumptions mentioned above,
$$\lim\li_{\dt \ra 0}\E \max\li_{t\in
[0,T]}|\widetilde{\boldq}_t^\dt-\widehat{\boldq}_t|_{\R^d}=0  \ , $$
where $\widehat{\boldq}_t$ is the solution of the problem}

$$\dot{\widehat{\boldq}}_t=\dfrac{\boldb(\widehat{\boldq}_t)}{\lb(\widehat{\boldq}_t)}+
\dfrac{\sm(\widehat{\boldq}_t)}{\lb(\widehat{\boldq}_t)}\circ
\dot{\boldW}_t \ , \ \widehat{\boldq}_0=\boldq\in \R^d \ .
\eqno(3.9)$$

\

\section{One dimensional case}

In the case of one space variable, Smoluchowski-Kramers
approximation leads to an one-dimensional diffusion process $q_t$
which is defined by the following stochastic differential equation
written in the It\^{o} form:

$$\dot{q}_t=\dfrac{b(q_t)}{\lb(q_t)}-\dfrac{\lb'(q_t)}{2\lb^3(q_t)}
+\dfrac{1}{\lb(q_t)}\dot{W}_t \ , \ q_0=q\in \R^1 \ . \eqno(4.1)$$

Put

$$\begin{array}{l}
\play{u(q)=\int_0^q \lb(x) \exp \left(-2\int_0^x
b(y)\lb(y)dy\right)dx \ ,}
\\
\play{v(q)=2\int_0^q \lb(x) \exp\left(2\int_0^x
b(y)\lb(y)dy\right)dx \ .}
\end{array} \eqno(4.2)$$

Since $\lb(x)>0$, $u(q)$ and $v(q)$ are strictly increasing
functions. Following [4] we introduce an operator $D_vD_u$, where
$D_u$ means the differentiation with respect to the monotone
function $u(q)$: $\play{D_u f(q)=\lim\li_{h\ra
0}\dfrac{f(x+h)-f(x)}{u(x+h)-u(x)}}$; the operator $D_v$ is defined
in a similar way. One can check that $D_vD_u$ is the generator of
the diffusion process $q_t$ defined by (4.1).

Suppose now that the friction coefficient $\lb(q)=\lb_\ve(q)$
depends on a parameter $\ve > 0$. We assume that, for each $\ve \in
(0,1]$, $\lb_\ve(q)$ has a bounded continuous derivative
$\lb_\ve'(q)$, and $0< \underline{\lb}\leq \lb_\ve(q) \leq
\overline{\lb} <\infty$. Let $u_\ve(q)$ and $v_\ve(q)$ be the
functions defined by (4.2) when $\lb(q)$ is replaced by
$\lb_\ve(q)$.

Consider the stochastic process $q_t^{\mu,\dt,\ve}$ in $\R^1$
defined by the equation

$$\mu \ddot{q}_t^{\mu,\dt,\ve}=b(q_t^{\mu,\dt,\ve})-
\lb^\ve(q_t^{\mu,\dt,\ve})q_t^{\mu,\dt,\ve}+\dot{W}_t^\dt \ , \
q_0^{\mu,\dt,\ve}=q \ , \ \dot{q}_0^{\mu,\dt,\ve}=p \ . \eqno(4.3)$$
where $\dot{W}_t^\dt$ is, as before, a "smoothed" white noise
converging to $\dot{W}_t$ as $\dt \da 0$.

\

\textbf{Theorem 4.1.} \textit{Assume that the function $\lb_\ve(q)$
converge weakly as $\ve \da 0$ on each finite interval $[\al,
\bt]\subset \R^1$ to a function $\overline{\lb}(q)$ (maybe,
discontinuous). Then processes $q_t^{\mu,\dt,\ve}$ converge weakly
on each finite time interval to the diffusion process
$\overline{q}_t$ governed by the generator
$D_{\overline{v}}D_{\overline{u}}$ (where $\overline{u}(q)$ and
$\overline{v}(q)$ defined by (4.2) with $\lb=\overline{\lb}(q)$) as,
first $\mu \da 0$, then $\dt \da 0$, and then $\ve \da 0$.}

\

\textbf{Proof.} According to Section 3, processes
$q_t^{\mu,\dt,\ve}$ converge weakly as first $\mu \da 0$ and then
$\dt \da 0$ to the process $\widehat{q}_t^\dt$ which solves equation
(4.1) with $\lb(q)=\lb^\ve(q)$. It follows from our assumptions that
functions $u_\ve(q)$ and $v_\ve(q)$ converge as $\ve \da 0$ to
functions $\overline{u}(q)$ and $\overline{v}(q)$ respectively for
each $q\in \R^1$. The functions $\overline{u}(q)$ and
$\overline{v}(q)$ are continuous and strictly increasing. Therefore
([4]) a diffusion process $\overline{q}_t$ exists governed by
$D_{\overline{v}}D_{\overline{u}}$. As shown in [8], convergence of
$u_\ve(q)$ and $v_\ve(q)$ as $\ve \da 0$ to $\overline{u}(q)$ and
$\overline{v}(q)$ respectively implies weak convergence of processes
$q_t^\ve$ to the process corresponding to
$D_{\overline{v}}D_{\overline{u}}$ as $\ve \da 0$. $\square$

\

\textbf{Theorem 4.2. } \textit{Let
$\lb_\ve(q)=\widetilde{\lb}\left(\dfrac{q}{\ve}\right)$. Assume that
one of the following conditions is satisfied:}

\textit{1. $\widetilde{\lb}(q)$ is a continuously differentiable
positive 1-periodic function; }

\textit{2. $\widetilde{\lb}(q)$ is an ergodic stationary process
(independent of the process $W_t$ in (4.3)) with continuously
differentiable trajectories and $0<\lb_{-}\leq
\widetilde{\lb}(q)\leq \lb_{+}<\infty$ for some constants $\lb_{-}$,
$\lb_{+}$. }

\textit{Put $\overline{\lb}=\play{\int_0^1\widetilde{\lb}(q)dq}$ if
condition 1 is satisfied, and $\overline{\lb}=\E\widetilde{\lb}(q)$
is condition 2 is satisfied. }

\textit{Then the process $q_t^{\mu,\dt,\ve}$ defined by (5.3)
converge weakly when first $\mu \da 0$ and then $\ve \da 0$ to the
process $\overline{q}_t$ defined by the equation
$$\overline{q}_t=\dfrac{1}{\overline{\lb}}b(\overline{q}_t)+
\dfrac{1}{\overline{\lb}}\dot{W}_t \ , \ \overline{q}_0=q \ .$$ }

\

\textbf{Proof} of this theorem follows from Theorem 4.1 since each
of conditions 1 and 2 implies conditions of Theorem 4.1 and
$\overline{\lb}(q)=\overline{\lb}$.  $\square$

\

Assume now that $\lb_{\ve}(q)$ is a bounded and separated from zero
uniformly in $\ve\in (0,1]$ positive function such that
$\lim\li_{\ve\da 0}\lb_\ve(q)=\lb_1$ for $q<0$, and $\lim\li_{\ve
\da 0}\lb_\ve(q)=\lb_2$ for $q>0$. Assume that $\lb_\ve(q)$ is
continuously differentiable for each $\ve > 0$. Let
$\widehat{\lb}(q)$ be the step function equal to $\lb_1$ for $q\leq
0$ and to $\lb_2$ for $q>0$. Let functions $\widehat{u}(q)$ and
$\widehat{v}(q)$ be defined by formula (4.2) with
$\lb(q)=\widehat{\lb}(q)$; $\widehat{u}(q)$ and $\widehat{v}(q)$ are
continuous strictly increasing functions. Denote by $\widehat{q}_t$
the diffusion process in $\R^1$ governed by the generator
$A=D_{\widehat{v}}D_{\widehat{u}}$. The process $\widehat{q}_t$
behaves as $\dfrac{1}{\lb_1}W_t$ on the negative part of axis $q$
and as $\dfrac{1}{\lb_2}W_t$ on the positive part. Its behavior at
$q=0$ is defined by the domain of definition $\mathfrak{D}_A$ of the
generator $A$: a continuous bounded function $f(q)$, $q\in \R^1$,
twice continuously differentiable at $q\in \{\R^1 \setminus
\{q=0\}\}$ belongs to $\mathfrak D_A$ if and only if left and right
derivatives at $q=0$, $f_{-}'(0)$ and $f_{+}'(0)$ respectively,
satisfy the equality
$\dfrac{1}{\lb_1}f'_{-}(0)=\dfrac{1}{\lb_2}f'_{+}(0)$ and $Af(q)$ is
continuous.

It is easy to see that functions $u_\ve(q)$ and $v_\ve(q)$ defined
by (4.2) with $\lb(q)=\lb_\ve(q)$ converge as $\ve \da 0$ to
$\widehat{u}(q)$ and $\widehat{v}(q)$ respectively for each $q\in
\R^1$. This implies the following result.

\

\textbf{Theorem 4.3.} \textit{Let the friction coefficient
$\lb_\ve(q)$ satisfies the conditions mentioned above. Then the
stochastic process $q_t^{\mu,\dt,\ve}$ defined by (4.3) converges
weakly to the diffusion process $\widehat{q}_t$ in $\R^1$ governed
by $A=D_{\widehat{v}}D_{\widehat{u}}$ as first $\mu \da 0$, then
$\dt \da 0$, and then $\ve \da 0$.}

\

This means, roughly speaking, that, if the friction coefficient is
close to the step-function $\widehat{\lb}(q)$, then process
$q_t^\mu$, for $0<\mu<<1$, can be approximated by the diffusion
process $\widehat{q}_t$.

\

\section{Multidimensional case}

In this section we consider the problem of fast oscillating periodic
environment in multidimensional case. We consider the system

$$\mu \ddot{\boldq}^{\mu,\dt,\ve}_t=\boldb\left(
\dfrac{\boldq_t^{\mu,\dt,\ve}}{\ve}\right)
-\lb\left(\dfrac{\boldq_t^{\mu,\dt,\ve}}
{\ve}\right)\dot{\boldq}_t^{\mu,\dt,\ve} +\dot{\boldW}_t^\dt \ , \
\boldq_0^{\mu,\dt,\ve}=\boldq\in \R^d
 \ , \ \dot{\boldq}_0^{\mu,\dt,\ve}=\boldp\in \R^d \ . \eqno(5.1)$$

Here as in Section 3 the process $\boldW_t^\dt$ is the approximation
of the Wiener process in $\R^d$. We make the same assumptions about
the functions $\lb(\bullet)$ and $\boldb(\bullet)$ as in Section 2.
In addition we assume that the functions $\lb(\bullet)$ and
$\boldb(\bullet)$ are 1-periodic, i.e.
$\lb(\boldx+\mathbf{e}_k)=\lb(\boldx)$ and
$\boldb(\boldx+\mathbf{e}_k)=\boldb(\boldx)$ for $\boldx\in \R^d$
and $\mathbf{e}_k=(0,0,...,1 (k\text{-th coordinate}),...,0)$,
$1\leq k \leq d$. Under this assumption our system (5.1) could be
regarded as a system on the $d$-torus $\torus^d=\R^d/\Z^d$. Fix
$\ve>0$, we can proceed as in Section 3 to see that first as $\mu
\da 0$ then as $\dt \da 0$ the process $\boldq_t^{\mu,\dt,\ve}$
converges in probability to the process $\boldq_t^{\ve}$ subjected
to

$$\dot{\boldq}_t^\ve=\dfrac{\boldb\left(
\dfrac{\boldq_t^\ve}{\ve}\right)}{\lb\left(\dfrac{\boldq_t^\ve}{\ve}\right)}+
\dfrac{1}{\lb\left(\dfrac{\boldq_t^\ve}{\ve}\right)}\circ
\dot{\boldW}_t \ , \boldq_0^\ve=\boldq\in \R^d \ .
$$

The above equation, written in the form of It\^{o} integral, will be

$$\dot{\boldq}_t^\ve=\dfrac{\boldb\left(
\dfrac{\boldq_t^\ve}{\ve}\right)}{\lb\left(\dfrac{\boldq_t^\ve}{\ve}\right)}-
\dfrac{1}{2\ve} \dfrac{\grad
\lb\left(\dfrac{\boldq_t^\ve}{\ve}\right)}{\lb^3\left(\dfrac{\boldq_t^\ve}{\ve}\right)}+
\dfrac{1}{\lb\left(\dfrac{\boldq_t^\ve}{\ve}\right)}\dot{\boldW}_t \
, \boldq_0^\ve=\boldq\in \R^d \ . \eqno(5.2)
$$

The generator corresponding to (5.2) is the second order
differential operator

$$L^\ve u(\boldx)=\left(\dfrac{\boldb\left(\dfrac{\boldx}{\ve}\right)}
{\lb\left(\dfrac{\boldx}{\ve}\right)}- \dfrac{1}{2\ve}\dfrac{\grad
\lb
\left(\dfrac{\boldx}{\ve}\right)}{\lb^3\left(\dfrac{\boldx}{\ve}\right)}\right)\cdot
\grad u(\boldx)
+\dfrac{1}{2}\dfrac{1}{\lb^2\left(\dfrac{\boldx}{\ve}\right)}\Dt
u(\boldx) \ . \eqno(5.3)$$

Our goal is to study the homogenization properties of (5.3) for
general multidimensional case. Homogenization problems are
considered by many authors, see, e.g., [7], [15], [14], [13], [12].
However, we provide here an elementary probabilistic way of doing
this. Our method follows [7] and [6](pp. 104-106).

Let us first make a change of variable $\dfrac{\boldq}{\ve}=\boldy$
and $\dfrac{t}{\ve^2}=s$. The process
$\boldy_s^\ve=\dfrac{1}{\ve}\boldq_t^\ve$ corresponds to the
generator
$$A^\ve=\dfrac{1}{2\lb^2(\boldy)}\Dt_{\boldy}-\dfrac{\grad \lb(\boldy)}{2\lb^3(\boldy)}
\cdot \grad_{\boldy} +\ve\dfrac{\boldb(\boldy)}{\lb(\boldy)}\cdot
\grad_{\boldy} \ .$$

We regard $\boldy_s^\ve$ as a process on $\torus^d$. Then we have
the bound
$$\left|\E_{\boldq/\ve}f(\boldy_s^\ve)-\int_{\torus^d}f(\boldx)
\mu^\ve(\boldx)d\boldx\right|<Ke^{-as} \ .$$

Here $K>0$ and $a>0$ are independent of $\ve$ for small $\ve$. The
function $f$ is bounded and measurable. The function
$\mu^\ve(\boldx)$ is the density of the unique invariant measure of
$\boldy_s^\ve$ on $\torus^d$ and
$\play{\int_{\torus^d}\mu^\ve(\boldx)d\boldx=1}$. We have
$$\lim\li_{\ve\da 0}\mu^\ve(\boldx)=\mu(\boldx) \ ,
\ \lim\li_{\ve \da
0}\int_{\torus^d}f(\boldx)\mu^\ve(\boldx)d\boldx=\int_{\torus^d}f(\boldx)\mu(\boldx)d\boldx$$
for $f\in C(\torus^d)$ and $\mu(\boldx)$ the unique invariant
measure for the process with generator $A^0$ on $\torus^d$ and
$\play{\int_{\torus^d}\mu(\boldx)d\boldx=1}$. Combining these
estimates we have, that for any $n$, for any $t\geq \dt >0$, there
exist $\ve_0(n,\dt)>0$ such that for any $0<\ve<\ve_0(n,\dt)$, we
have

$$\left|\E_{\boldq}f\left(\dfrac{\boldq_t^\ve}{\ve}\right)-
\int_{\torus^d}f(\boldx)\mu(\boldx)d\boldx\right|<\dfrac{1}{n} \ .$$
This implies that for any $f\in C(\torus^d)$,

$$\lim\li_{\ve\da 0}\sup\li_{t \geq \dt}
\left|\E_{\boldq}f\left(\dfrac{\boldq_t^\ve}{\ve}\right)-
\int_{\torus^d}f(\boldx)\mu(\boldx)d\boldx\right|=0 \ .$$

Finally we calculate the density $\mu(\boldx)$. Since
$$A^0=\dfrac{1}{2\lb^2(\boldy)}\Dt_{\boldy}-\dfrac{\grad \lb(\boldy)}{2\lb^3(\boldy)}\cdot \grad_{\boldy}=
\dfrac{1}{2\lb^2(\boldy)}(\Dt_{\boldy}-\grad(\ln \lb(\boldy))\cdot
\grad_{\boldy}) \ ,$$ we see that $\mu(\boldx)=C\lb(\boldx)$ with
$C=\left(\play{\int_{\torus^d}\lb(\boldx)d\boldx}\right)^{-1}$ and
we have the following result:

\

\textbf{Lemma 5.1.} \textit{For any $f\in C(\torus^d)$, we have
$$\lim\li_{\ve\da 0}\sup\li_{t\geq \dt}\left|\E_{\boldq}
f\left(\dfrac{\boldq_t^\ve}{\ve}\right)-\dfrac{\play{\int_{\torus^d}f(\boldx)\lb(\boldx)d\boldx}}
{\play{\int_{\torus^d}\lb(\boldx)d\boldx}}\right| = 0 \ .
\eqno(5.4)$$ }

\

\textbf{Corollary.} \textit{ For any bounded continuous function
$f(\boldx)$ on $\torus^d$, $\boldq\in \torus^d$ we have}

$$\E_{\boldq}\left[\int_0^t f\left(\dfrac{\boldq_s^\ve}{\ve}\right)ds-
\dfrac{t\play{\int_{\torus^d}f(\boldx)\lb(\boldx)d\boldx}}
{\play{\int_{\torus^d}\lb(\boldx)d\boldx}}\right]^2\ra 0$$
\textit{as $\ve \da 0$, for $0<t<\infty$.}

\

The \textit{proof} of this corollary follows the same proof of the
corollary after Lemma 1 in [7].

\

Now let us consider auxiliary functions $N_k(\boldy)$, $k=1,...,d$,
which are the periodic bounded solutions (i.e., on $\torus^d$) of
the equations

$$\dfrac{1}{2\lb^2(\boldy)}\Dt_{\boldy} N_k(\boldy)-\dfrac{\grad_{\boldy} \lb(\boldy)}{2\lb^3(\boldy)}\cdot
\grad_{\boldy}
N_k(\boldy)=A^0(N_k(\boldy))=\dfrac{1}{2\lb^3(\boldy)}\play{\dfrac{\pt
\lb}{\pt y_k}(\boldy)} \ , \ \boldy\in \torus^d  \ . \eqno(5.5)$$

The solvability of this equation comes from the fact that $(A^0)^*
\lb(\boldy)=0$ and
$\play{\int_{\torus^d}\dfrac{1}{2\lb^3(\boldy)}\dfrac{\pt \lb}{\pt
y_k}(\boldy)\lb(\boldy)d\boldy}=0$. The boundedness of solution
comes from our assumptions on the function $\lb(\bullet)$. Now we
apply It\^{o}'s formula:

$$\begin{array}{l}
\ve N_k\left(\dfrac{\boldq_t^\ve}{\ve}\right)-\ve
N_k\left(\dfrac{\boldq}{\ve}\right)
\\
\play{= \ve\left[\int_0^t \grad
N_k\left(\dfrac{\boldq_s^\ve}{\ve}\right)\cdot\dfrac{1}{\ve}
\left(\dfrac{\boldb}{\lb}\left(\dfrac{\boldq_s^\ve}{\ve}
\right)-\dfrac{1}{2\ve}\dfrac{\grad
\lb}{\lb^3}\left(\dfrac{\boldq_s^\ve}{\ve}\right)+\dfrac{\dot{\boldW}_s}
{\lb\left(\dfrac{\boldq_s^\ve}{\ve}\right)}\right)ds+\right.}
\\
\play{ \ \ \ \ \ \ \ \ \ \ \ \ \ \ \left. + \dfrac{1}{2}\int_0^t \Dt
N_k\left(\dfrac{\boldq_s^\ve}{\ve}\right)\dfrac{1}{\ve^2}\dfrac{1}{\lb^2
\left(\dfrac{\boldq_s^\ve}{\ve}\right)}ds\right]}
\\
\play{=\int_0^t \grad N_k\left(\dfrac{\boldq_s^\ve}{\ve}\right)\cdot
\left(\dfrac{\boldb}{\lb}\left(\dfrac{\boldq_s^\ve}{\ve}
\right)+\dfrac{\dot{\boldW}_s}
{\lb\left(\dfrac{\boldq_s^\ve}{\ve}\right)}\right)ds+\dfrac{1}{2\ve}\int_0^t
\dfrac{\dfrac{\pt \lb}{\pt
y_k}\left(\dfrac{\boldq_s^\ve}{\ve}\right)}{\lb^3
\left(\dfrac{\boldq_s^\ve}{\ve}\right)}ds} \ .
\end{array} \eqno(5.6)$$

Let $\boldN(\boldy)=(N_1(\boldy), ... , N_d(\boldy))$. Using (5.5)
we have

$$\begin{array}{l}
\play{\boldq_t^\ve-\boldq=\int_0^t\left(\dfrac{\boldb}{\lb}\left(\dfrac{\boldq_s^\ve}{\ve}\right)+
\dfrac{\dot{\boldW}_s}{\lb\left(\dfrac{\boldq_s^\ve}{\ve}\right)}\right)ds
+\int_0^t (D\boldN)\left(\dfrac{\boldq_s^\ve}{\ve}\right)
\left(\dfrac{\boldb}{\lb}\left(\dfrac{\boldq_s^\ve}{\ve}
\right)+\dfrac{\dot{\boldW}_s}
{\lb\left(\dfrac{\boldq_s^\ve}{\ve}\right)}\right)ds + }
\\
\play{ \ \ \ \ \ \ \ \ \ \ \ \ \ \ \ \ \  - \ve
\left(\boldN\left(\dfrac{\boldq_t^\ve}{\ve}\right)-
\boldN\left(\dfrac{\boldq}{\ve}\right)\right)} \ ;
\end{array}$$

Here $(D\boldN)(\boldy)=\left(\dfrac{\pt N_i}{\pt y_j}\right)_{1\leq
i,j\leq d}$,
 $\boldy=(y_1,...,y_d)\in \torus^d$.

Therefore using the corollary after Lemma 1, we see that
$\boldq_t^\ve$ converges weakly to a process $\boldq_t$,
$\boldq_0=\boldq\in \R^d$ governed by the operator

$$\overline{L}=\dfrac{1}{2}\sum\li_{i,j=1}^d \bar{a}_{ij}\dfrac{\pt^2}{\pt y_i\pt y_j}+\sum\li_{i=1}^d
\bar{b}_i \dfrac{\pt}{\pt y_i} \ , \eqno(5.7)$$ with coefficients

$$\bar{a}_{ij}=\int_{\torus^d}\left(\dfrac{\grad N_i(\boldy)\cdot\grad N_j(\boldy)}{\lb(\boldy)}
+ \dfrac{1}{\lb(\boldy)}\left(\dfrac{\pt N_j}{\pt y_i}(\boldy) +
\dfrac{\pt N_i}{\pt y_j}(\boldy)\right)
+\dt_{ij}\dfrac{1}{\lb(\boldy)}\right)d\boldy \left/
\left(\int_{\torus^d}\lb(\boldy)d\boldy\right) \right. \ , $$

$$\bar{b}_i=\dfrac{\play{\int_{\torus^d}b_i(\boldy)d\boldy}}{\play{\int_{\torus^d}\lb(\boldy)d\boldy}}
+ \sum\li_{k=1}^d \dfrac{\play{\int_{\torus^d}b_k(\boldy)\dfrac{\pt
N_i}{\pt
y_k}(\boldy)d\boldy}}{\play{\int_{\torus^d}\lb(\boldy)d\boldy}} \ .
\eqno(5.8)$$

Here $\dt_{ij}=1$ if $i=j$, and $\dt_{ij}=0$ otherwise.

We could simplify the expression for $\bar{a}_{ij}$: using (5.5) we
get

$$\begin{array}{l}
\play{\bar{a}_{ij}=\int_{\torus^d}\left(\dfrac{\grad
N_i(\boldy)\cdot\grad N_j(\boldy)}{\lb(\boldy)} +
\dfrac{1}{\lb(\boldy)}\left(\dfrac{\pt N_j}{\pt y_i}(\boldy) +
\dfrac{\pt N_i}{\pt y_j}(\boldy)\right)
+\dt_{ij}\dfrac{1}{\lb(\boldy)}\right)d\boldy\left/
\left(\int_{\torus^d}\lb(\boldy)d\boldy\right) \right.}
\\
\play{=\int_{\torus^d}\left(\text{div}\left(\dfrac{N_i(\boldy)}{\lb(\boldy)}\grad
N_j(\boldy)\right) - \dfrac{N_i(\boldy)}{\lb(\boldy)}\Dt N_j
(\boldy) - N_i(\boldy) \grad N_j(\boldy)
 \cdot \grad \left(\dfrac{1}{\lb(\boldy)}\right) + \right.}
\\
\play{\left. \ \ \ \ \ \ \ \ \ \ \ \ \ \ \ +
\dfrac{1}{\lb(\boldy)}\left(\dfrac{\pt N_j}{\pt y_i}(\boldy) +
\dfrac{\pt N_i}{\pt y_j}(\boldy)\right)
+\dt_{ij}\dfrac{1}{\lb(\boldy)}\right)d\boldy\left/
\left(\int_{\torus^d}\lb(\boldy)d\boldy\right) \right.}
\\
\play{=\int_{\torus^d}\left(\dfrac{\pt}{\pt
y_j}\left(N_i(\boldy)\dfrac{1}{\lb(\boldy)}\right) -
\dfrac{1}{\lb(\boldy)}\dfrac{\pt N_i}{\pt y_j}(\boldy)  \ + \right.}
\\
\play{\left. \ \ \ \ \ \ \ \ \ \ \ \ \ \ \ +
\dfrac{1}{\lb(\boldy)}\left(\dfrac{\pt N_j}{\pt y_i}(\boldy) +
\dfrac{\pt N_i}{\pt y_j}(\boldy)\right)
+\dt_{ij}\dfrac{1}{\lb(\boldy)}\right)d\boldy\left/
\left(\int_{\torus^d}\lb(\boldy)d\boldy\right) \right.}
\\
\play{=\dfrac{\play{\int_{\torus^d}\dfrac{\pt N_j}{\pt
y_i}(\boldy)\dfrac{1}{\lb(\boldy)}d\boldy}}
{\play{\int_{\torus^d}\lb(\boldy)d\boldy}}+\dt_{ij}\dfrac{\play{\int_{\torus^d}\dfrac{1}{\lb(\boldy)}d\boldy}}
{\play{\int_{\torus^d}\lb(\boldy)d\boldy}} \ .}
\end{array} \eqno(5.9)$$

So we have

\

\textbf{Theorem 5.1.} \textit{As $\ve \da 0$, the process
$\boldq_t^\ve$ converges weakly to a process $\boldq_t$,
$\boldq_0=\boldq\in \R^d$ governed by the operator (5.7) with
coefficients given by (5.8) and (5.9). }

\

This Theorem implies a homogenization result for the process
$\boldq_t^{\mu,\dt,\ve}$ defined by equation (5.1).

\

\section{Remarks and Generalizations}

\textbf{6.1. Small mass - large friction asymptotics.}

Let the friction coefficient in (1.1) be
$\lb^\ve(\boldq)=\ve^{-1}\lb(\boldq)$, $\boldq\in \R^n$, $0<\ve<<1$.
As it follows from Theorem 3.1, the Smoluchowski-Kramers
approximation in this case has the form:

$$\dot{\overline{\boldq}}_t^\ve=\dfrac{\ve\boldb(\overline{\boldq}_t^\ve)}
{\lb(\overline{\boldq}_t^\ve)}-\dfrac{\ve^2 \grad
\lb(\overline{\boldq}_t^\ve)}{2\lb^3(\overline{\boldq}_t^\ve)} +
\dfrac{\ve}{\lb(\overline{\boldq}_t^\ve)}\dot{\boldW}_t \ , \
\overline{\boldq}_0^\ve=\boldq \ .$$

Put $\widetilde{\boldq}_t^\ve=\overline{\boldq}_{t/\ve}^\ve$. Then
$\widetilde{\boldq}_t^\ve$ satisfies the equation

$$\dot{\widetilde{\boldq}}_t^\ve=\dfrac{\boldb(\widetilde{\boldq}_t^\ve)}{\lb(\widetilde{\boldq}_t^\ve)}-
\dfrac{\ve\grad
\lb(\widetilde{\boldq}_t^\ve)}{2\lb^3(\widetilde{\boldq}_t^\ve)}
+\dfrac{\sqrt{\ve}}{\lb(\widetilde{\boldq}_t^\ve)}\dot{\widetilde{\boldW}}_t
\ , \ \widetilde{\boldq}_0^\ve=\boldq \ , \eqno(6.1)$$ where
$\widetilde{\boldW}_t$ is a Wiener process.

Assume that the vector field $\boldb(\boldq)$, $\boldq\in \R^n$, has
a finite number of compact attractors $K_1,...,K_l$. Let, for
brevity, each $K_i$ be an asymptotically stable equilibrium, and
each point of $\R^n$, besides a separatrix set $\cE\subset \R^n$, is
attracted to one of these equilibriums. The separatrix set $\cE$ is
assumed to have dimension less than $n$. Then, if
$\widetilde{\boldq}_0^\ve=\boldq\not \in \cE$,
$\widetilde{\boldq}_t^\ve$ first comes to a small neighborhood of a
stable equilibrium $K_i$, $i=i(\boldq)$, with the probability close
to $1$ as $\ve \da 0$ and spends in this neighborhood a long time.
Then, because of the large deviations, the trajectory will switch to
the neighborhood of another attractor, then to another, and so on.
We see from (6.1), that the long-time behavior of the system with a
large friction is similar to the behavior of a system with a small
noise. Applying the results of [9], Chapters 4, 6, we see that, for
$0<\ve<<1$, the sequence of transitions between the attractors, the
main term of transition time logarithmic asymptotics, and the most
probable transition paths are not random for generic systems. These
characteristics of the long-time behavior are defined by a function
$V(\boldx,\boldy)$:

$$V(\boldx,\boldy)=\inf\left\{\dfrac{1}{2}\int_0^T
|\lb(\boldphi_s)\dot{\boldphi}_s-\boldb(\boldphi_s)|^2ds :
\boldphi_0=\boldx \ , \ \boldphi_T=\boldy , T\geq 0  \right\}
 \ , \ \boldx,\boldy\in \R^n \ ;$$
 and by the extremals of this variational problem.

Let now the dynamical system $\dot{\boldq}_t=\boldb(\boldq_t)$ has a
first integral. Let, say, $n=2$ and
$\boldb(\boldq)=\overline{\grad}H(\boldq)$ for some smooth generic
function $H(\boldq)$, $\boldq\in \R^2$, such that
$\lim\li_{|\boldq|\ra \infty}H(\boldq)=\infty$; in this case
$H(\boldq_t)\equiv H(\boldq_0)$.

Assume again that the friction is strong:
$\lb^\ve(\boldq)=\ve^{-1}\lb(\boldq)$. Make a time change
$\widehat{\boldq}_t^\ve=\boldq_{t/\ve^2}$. Then

$$\dot{\widehat{\boldq}}_t^\ve=\dfrac{1}{\ve^2}\grad H(\widehat{\boldq}_t^\ve)
-\dfrac{\grad
\lb(\widehat{\boldq}_t^\ve)}{2\lb^3(\widehat{\boldq}_t^\ve)}
+\dfrac{1}{\lb(\widehat{\boldq}_t^\ve)}\dot{\boldW}_t \ , \
\widehat{\boldq}_0^\ve=\boldq\in \R^2 \ .$$

Identify points of each connected component of every level set of
$H(\boldq)$. The set obtained after such an identification if
homeomorphic in the natural topology to a graph $\Gm$. Let $Y:
\R^2\ra \Gm$ be the identification mapping. Then the long-time
evolution of the system can be characterized by the stochastic
process $\mathcal{Y}_t^\ve=Y(\widehat{\boldq}_t^\ve)$ on $\Gm$. The
process $\mathcal{Y}_t^\ve$, in general, is not Markovian. But
$\mathcal{Y}_t^\ve$ converges weakly in the space of continuous
functions $\varphi: [0,T]\ra \Gm$ as $\ve \da 0$ to a diffusion
process on the graph $\Gm$ ([9], Ch.8). This limiting process is
defined by a family of second-order differential operators, one on
each edge of $\Gm$, and by the gluing conditions at the vertices.
Following [9], one can evaluate these operators and the gluing
conditions.

\

\textbf{6.2. Fast oscillating random friction in multidimensional
case.}

Let us consider the case of fast oscillating in the space variable,
random friction, in dimension $d\geq 2$. Let $(\Om, \cF, \Prob)$ be
a probability space. Let $\lb(\boldx, \om)$, $\om \in \Om$ be a
random field in $\R^d$ with the following properties:

(i) For any fixed $\om\in \Om$ and $\boldx \in \R^d$, the function
$\infty>\Lb\geq\lb(\boldx,\om)\geq \lb_0>0$.

(ii) For every $\boldx\in \R^d$ the random variable $\lb(\boldx,
\om)$ is independent of the Wiener process $\boldW_t$.

(iii) The random field $\lb(\boldx,\om)$ has the form
$\lb(\boldx,\om)=\lb(T(\boldx)\om)$ where $T(\boldx): \Om\ra \Om$ is
a $d$-dimensional dynamical system which preserves the measure
$\Prob$ and is ergodic with respect to $\Prob$.

Let us now consider an analogue of (5.1):

$$\mu\ddot{\boldq}_t^{\mu,\dt,\ve}=-\lb
\left(\dfrac{\boldq_t^{\mu,\dt,\ve}}{\ve},\om\right)\dot{\boldq}_t^{\mu,\dt,\ve}
+\dot{\boldW}_t^\dt \ , \ \boldq_0^{\mu,\dt,\ve}=\boldq\in \R^d \ ,
\ \dot{\boldq}_0^{\mu,\ve,\dt}=\boldp\in \R^d \ . \eqno(6.2)$$

For each fixed $\om \in \Om$, as we proved in Section 3, we have
that $\boldq_t^{\mu,\dt,\ve}(\om)$ converges weakly to a process
$\boldq_t^\ve(\om)$ as first $\mu \da 0$ and then $\dt \da 0$. The
process $\boldq_t^\ve$ is subject to

$$\dot{\boldq}_t^\ve=-\dfrac{1}{2\ve}\dfrac{\grad \lb \left(\dfrac{\boldq_t^\ve}{\ve},\om\right)}
{\lb^3\left(\dfrac{\boldq_t^\ve}{\ve},\om\right)} +
\dfrac{1}{\lb\left(\dfrac{\boldq_t^\ve}{\ve},\om\right)}\dot{\boldW}_t
\ , \ \boldq_0^\ve=\boldq\in \R^d \ . \eqno(6.3)$$

We conjecture that as $\ve\da 0$, $\boldq_t^\ve$ converges weakly to
a process $\boldq_t$, $\boldq_0=\boldq\in \R^d$ subject to the
operator $\overline{\overline{L}}=\dfrac{1}{2}\sum\li_{i,j=1}^d
\overline{\overline{a}}_{ij}\dfrac{\pt^2}{\pt x_i \pt x_j}$ with
effective diffusivity
$$\overline{\overline{a}}_{ij}=\E\left[\dfrac{\play{\int_{\torus^d}\dfrac{\pt
N_j}{\pt y_i}(\boldx,\om)\dfrac{1}{\lb(\boldx,\om)}d\boldx}}
{\play{\int_{\torus^d}\lb(\boldx,\om)d\boldx}}+\dt_{ij}\dfrac{\play{\int_{\torus^d}\dfrac{1}{\lb(\boldx,\om)}d\boldx}}
{\play{\int_{\torus^d}\lb(\boldx,\om)d\boldx}}\right] \ .
\eqno(6.4)$$

Here the functions $N_k(\boldx, \om)$ ($1\leq k \leq d$) shall
satisfy certain auxiliary problem. (We actually have a formulation
of this problem but we are not sure about its validity: we let
$N_k(\boldx,\om)$ be the solution of the equation

$$\E[(\grad_{\boldx}N_k(\boldx,\om)-\bolde_k)
\cdot \grad_{\boldx}\varphi(\boldx,\om)]=0 \ , $$ for all
$\varphi(\boldx,\om)$ smooth and compactly supported in $\boldx \in
\R^d$ and measurable with respect to $\om\in \Om$. The existence of
solutions to this problem is guaranteed by the Lax-Milgram lemma.)

However, we are not aware of the validity of this conjecture nor a
proof of it. We are also not sure about the correct reference of
such a problem. (We thank E.Kosygina for pointing out to us two
relevant papers [1] and [11].)

\

\textbf{6.3. Motion of charged particles in a magnetic field}

One can expect that using the regularization by smoothed white noise
one can get the Smoluchowski-Kramers approximation for the equation

$$\mu\ddot{\boldq}_t^\mu=\boldb(\boldq_t^\mu)-A(\boldq_t^\mu)\dot{\boldq}_t^\mu
+\sm \dot{\boldW}_t \ , \ \boldq_0^\mu=\boldq \ ,
\dot{\boldq}_0^\mu=\boldp \ ; \ \boldq, \boldp \in \R^n \ ,
\eqno(6.5)$$

Here $\sm>0$ and $A(\boldq)$ are a matrix-valued functions having
strictly positive eigenvalues for each $\boldq\in \R^n$. In
particular, if $A(\boldq)=\lb(\boldq)A$, where $\lb(\boldq)>0$ and
$A$ is a constant positive definite symmetric matrix, the problem
can be reduced to the case considered in Section 3 by an appropriate
linear change of variables.

If $A$ has a negative eigenvalue, the Smoluchowski-Kramers
approximation is not applicable. The case of $A$ with pure imaginary
eigenvalues is of interest since such equations describe the motion
of charged particles in a magnetic field. If $\lb=\text{const}$,
$n=2$ and $A=\begin{pmatrix}0&-1\\1&0\end{pmatrix}$, the problem was
considered in [3]. In this case the Smoluchowski-Kramers
approximation holds after a regularization. If
$\boldb(\boldq)=-\grad F(\boldq)$, $\boldq\in \R^2$, and
$A=\begin{pmatrix}0&-1\\1&0\end{pmatrix}$, one can show that the
regularization by the smoothed white noise leads to the equation

$$\boldq_t=\dfrac{1}{\lb(\boldq_t)}\overline{\grad} F(\boldq_t)-\sm^2
\dfrac{\overline{\grad}
\lb(\boldq)}{2\lb^3(\boldq)}+\dfrac{\sm}{\lb(\boldq)}\dot{\widetilde{\boldW}}_t
\ , \ \boldq_0=\boldq\in \R^2 \ . \eqno(6.6)$$

If the noise in (6.5) is small ($0<\sm<<1$) the motion described by
(6.6) has a fast and a slow components. Applying the results of [9],
[2], one can describe the limiting (as $\sm \da 0$) slow component
of $\widehat{\boldq}_t^\sm=\boldq_{t/\sm^2}$ as a diffusion process
on the graph corresponding to the potential $F(\boldq)$ (the graph
is homeomorphic to the set of connected components of the level sets
of $F(\boldq)$ provided with the natural topology).

\

\textbf{Acknowledgements}: This work is supported in part by NSF
Grants DMS-0803287 and DMS-0854982.

\

\section*{References}

[1]  G.Allaire, R.Orive, Homogenization of periodic non self-adjoint
problems with large drift and potential, \textit{Esaim-control
Optimisation and Calculus of Variations} , vol. \textbf{13}, no. 4,
pp. 735-749, 2007.

[2] M.Brin, M.Freidlin, On stochastic behavior of perturbed
Hamiltonian systems, \textit{Ergodic Theory and Dynamical Systems},
\textbf{20} (2000), No.1. 55-76.

[3] S.Cerrai, M.Freidlin, Small mass asymptotics of a charged
particle in a magnetic field and long-time influence of small
perturbations, \textit{Journal of Statistical Physics}, submitted.

[4] W.Feller, Generalized second order differential operators and
their lateral conditions, \textit{Ill. J. Math.}, \textbf{1} (1957)
pp.459-504.

[5] M.Freidlin, Some Remarks on the Smoluchowski-Kramers
Approximation, \textit{Journal of Statistical Physics},
\textbf{117}, No.314, pp.617-634, 2004.

[6] M.Freidlin, \textit{Functional Integration and Partial
Differential Equations}, Princeton University Press, 1985.

[7] M.Freidlin, Dirichlet's problem for an equation with periodic
coefficients depending on a small parameter, \textit{Theory of
Probability and Its Applications}, \textbf{9}, pp.133-139, 1964.

[8] M.Freidlin, A.Wentzell, Necessary and Sufficient Conditions for
Weak Convergence of One-Dimensional Markov Processes. \textit{The
Dynkin Festschrift: Markov Processes and their Applications} ,
Birkhauser, (1994) pp.95-109.

[9] M.Freidlin, A.Wentzell, \textit{Random perturbations of
dynamical systems}, Springer, 1998.

[10] N.Ikeda, S.Watanabe, \textit{Stochastic Differential Equations
and Diffusion Processes}, second edition, North-Holland Publishing
Company, 1989.

[11] M.Klepsyna, A.Piatnitski, Averaging of non-self adjoint
parabolic equations with random evolution, \textit{Institut National
de Recherche en Informatique et en Automatique}, June 2000.

[12] N.Krylov, G-convergence of elliptic operators in non-divergence
form, \textit{Math Zametki}, \textbf{37}, 4, pp. 522-527, April,
1985.

[13] G.Papanicolaou, S.R.S.Varadhan, Boundary value problems with
rapidly oscillating random coefficients, \textit{Seria Coll. Math.
Soc. J\'{a}nos Bolyai}, North Holland, Amsterdam, 1979.

[14] V.V.Zhikov, S.M.Kozlov, O.A.Oleinik, \textit{Homogenization of
differential operators and integral functionals}, Springer, 1994.

[15] V.V.Zhikov, S.M.Kozlov, O.A.Oleinik, and Kha T'en Ngoan,
Averaging and G-convergence of differential operators,
\textit{Russian Math Surveys}, \textbf{34}, 5, pp.65-133, 1979.

\end{document}